\documentclass[11pt]{article}
\usepackage[letterpaper,textwidth=6.3in,textheight=9in,centering]{geometry}
\usepackage[T1]{fontenc}
\usepackage{lmodern,microtype,amsmath,amssymb,amsthm,mathtools}
\usepackage{booktabs,tabularx,enumitem,xcolor}
\usepackage[hypertexnames=false,colorlinks=true,linkcolor=blue!40!black,citecolor=blue!40!black,urlcolor=blue!40!black]{hyperref}
\allowdisplaybreaks[2]
\setlist{nosep,leftmargin=1.7em}
\newtheorem{theorem}{Theorem}[section]
\newtheorem{proposition}[theorem]{Proposition}
\newtheorem{lemma}[theorem]{Lemma}
\newtheorem{corollary}[theorem]{Corollary}
\theoremstyle{definition}
\theoremstyle{remark}\newtheorem{remark}[theorem]{Remark}
\newcommand{\E}{\mathbb E}
\newcommand{\Pp}{\mathbb P}
\newcommand{\ET}{\mathbb E_T}
\newcommand{\PT}{\mathbb P_T}
\newcommand{\Var}{\operatorname{Var}}
\newcommand{\Cov}{\operatorname{Cov}}

\newcommand{\one}{\mathbf1}
\newcommand{\Rea}{\operatorname{Re}}
\newcommand{\Ima}{\operatorname{Im}}
\newcommand{\cC}{\mathcal C}

\newcommand{\cL}{\mathcal L}
\newcommand{\eps}{\varepsilon}

\newcommand{\Ac}{A_{\rm cert}}
\newcommand{\Cg}{C_{\rm geo}}
\title{Sharp tails for arguments of Dirichlet \(L\)-functions}
\author{Scott D. Hughes}
\date{}
\hypersetup{pdftitle={Sharp tails for arguments of Dirichlet L-functions},pdfauthor={Scott D. Hughes}}
\begin{document}
\maketitle
\vspace{-1em}
\begin{abstract}
Let \(t\) be uniform on \([T,2T]\), let \(\chi\) be a Dirichlet character
modulo \(q\), and set \(X=\Im\log L(1/2+it,\chi)\) and
\(Q=\sum_{p\le T,\,p\nmid q}1/p\). For fixed \(K,\theta_0>0\), uniformly
in \(q\le T^{\theta_0}\), \(\chi\), and \(0\le V\le KQ\), each signed
tail \(\Pr\{\pm X>V\}\) has order
\(\sqrt Q\,(V+\sqrt Q)^{-1}\exp(-V^2/Q)\).
In particular, for \(S(t)=\pi^{-1}\arg\zeta(1/2+it)\) and
\(0\le V\le K\log\log T\), each tail \(\Pr\{\pm S(t)>V\}\) has order
\(\sqrt{\log\log T}\,(V+\sqrt{\log\log T})^{-1}\exp(-\pi^2V^2/\log\log T)\).
For \(|v|\le KQ\), the probability that \(X\) lies in a fixed-width interval centered at
\(v\) has order \(Q^{-1/2}\exp(-v^2/Q)\) whenever the width is
sufficiently large. The upper bound holds for every fixed positive width.
We also obtain matching bounds of order \(\exp(k^2Q)\) for the means of
\(\exp(\pm2kX)\) for fixed \(k>0\).
The proof begins with a one-sided polynomial bound for each sign, using
Alp\"oge's zero-location theorem at the primitive inducing conductor.
A weighted comparison with a sum over primes controls departures from
the target interval, and a signed polynomial gives its lower bound.
\end{abstract}
\section{Results, analytic input, and proof structure}\label{sec:introduction}

Fix an integer \(q\ge1\) and a Dirichlet character \(\chi\bmod q\).
For \(q=1\), set \(L(s,\chi)=\zeta(s)\). Define
\[
 X_\chi(t)=\Im\log L(\tfrac12+it,\chi),\qquad
 \cL_T=\sum_{\substack{p\le T\\p\nmid q}}\frac1p,\qquad
 \ET f=\frac1T\int_T^{2T}f(t)\,dt,
\]
and write \(\PT\) for normalized Lebesgue measure on \([T,2T]\).
The logarithm is obtained by horizontal continuation from \(+\infty+it\),
with value tending to zero there. Values at zero ordinates can be
discarded. The primes dividing \(q\) are absent from the Euler product,
and \(\cL_T\) is the variance parameter that governs the tails:
\[
 \cL_T=\log\Bigl(\frac{\varphi(q)}q\log T\Bigr)+O(1),
\]
with an absolute \(O(1)\) whenever \(q\le T^{\theta_0}\) for a fixed
\(\theta_0>0\) and \(T\) is large (see the end of
Section~\ref{lw:sec:tails}).

The only information on zeros used in this paper is the zero-location
theorem of Alp\"oge~\cite{Alpoge} for \(\zeta\) and Dirichlet
\(L\)-functions, whose formal statements and
independent checks are recorded in~\cite{Sutherland}. It enters through
a single inequality, \textup{(OS)} below, and every result is first
proved from that inequality alone.

\begin{theorem}[Sharp tails for all Dirichlet characters]\label{cor:grh}
Fix \(K,\theta_0>0\). There are \(C<\infty\), \(c>0\) and \(T_0\),
depending only on \(K\) and \(\theta_0\), such that for \(T\ge T_0\),
every modulus \(1\le q\le T^{\theta_0}\), every Dirichlet character
\(\chi\bmod q\), \(0\le V\le K\cL_T\), and either \(\eps\in\{-1,1\}\),
\begin{equation}\label{eq:uniform}
 c\,\frac{\sqrt{\cL_T}}{V+\sqrt{\cL_T}}\,e^{-V^2/\cL_T}
 \le\PT\{\eps X_\chi>V\}\le
 C\,\frac{\sqrt{\cL_T}}{V+\sqrt{\cL_T}}\,e^{-V^2/\cL_T}.
\end{equation}
\end{theorem}

The common order is, up to constants, the tail of a Gaussian of
variance \(\cL_T/2\) throughout this range. In particular, uniformly
for \(\kappa\) in compact subsets of \((0,\infty)\),
\[
 \log\PT\{\eps X_\chi>\kappa\cL_T\}
 =-\kappa^2\cL_T-\tfrac12\log\cL_T+O(1).
\]
Replacing \(\cL_T\) by \(\log((\varphi(q)/q)\log T)\) in
\eqref{eq:uniform} changes only the constants. For a fixed modulus
\(\cL_T-\log\log T\) is bounded, and one may use \(\log\log T\) instead
(Corollary~\ref{thm:main}). When \(q\) grows,
\(\log\log T-\cL_T=\log(q/\varphi(q))+O(1)\) is unbounded along suitable
moduli (for example primorials), and at every fixed positive slope the
expression with \(\log\log T\) then differs from the true tail by an
unbounded factor.

For \(\zeta\) the theorem reads as follows in the usual normalization.

\begin{corollary}[The function \(S(t)\)]\label{cor:S}
Let \(S(t)=\pi^{-1}\Im\log\zeta(\tfrac12+it)\), with the branch fixed
above. For each fixed \(K>0\) there are \(C<\infty\), \(c>0\) and \(T_0\)
such that, for \(T\ge T_0\), \(0\le V\le K\log\log T\), and either sign,
\[
 c\,\frac{\sqrt{\log\log T}}{V+\sqrt{\log\log T}}\,
   e^{-\pi^2V^2/\log\log T}
 \le\PT\{\pm S(t)>V\}\le
 C\,\frac{\sqrt{\log\log T}}{V+\sqrt{\log\log T}}\,
   e^{-\pi^2V^2/\log\log T}.
\]
Moreover, for every real \(v\) with \(|v|\le K\log\log T\) and each fixed
\(h>0\),
\[
 \PT\{|S(t)-v|\le h\}\le C_h\,(\log\log T)^{-1/2}
   e^{-\pi^2v^2/\log\log T}\qquad(T\ge T_h),
\]
and there is \(h_1\), depending only on \(K\), such that for each fixed
\(h\ge h_1\) the reverse inequality holds with a positive constant in
place of \(C_h\).
\end{corollary}

\begin{proof}
Take \(q=1\). Then \(X_\chi=\pi S\), and Mertens' formula gives
\(\cL_T=\sum_{p\le T}p^{-1}=\log\log T+O(1)\). The event
\(\{\pm S>V\}\) is \(\{\pm X_\chi>\pi V\}\), and
\(\pi V\le(\pi K+1)\cL_T\) for large \(T\). Theorem~\ref{cor:grh} with
\(\pi K+1\) in place of \(K\) gives the order
\(\sqrt{\cL_T}\,(\pi V+\sqrt{\cL_T})^{-1}e^{-\pi^2V^2/\cL_T}\). Here
\[
 \pi^2V^2\Bigl|\frac1{\cL_T}-\frac1{\log\log T}\Bigr|
 =\frac{\pi^2V^2\,|\log\log T-\cL_T|}{\cL_T\log\log T}=O_K(1),
\]
and the prefactor is within a bounded factor of
\(\sqrt{\log\log T}\,(V+\sqrt{\log\log T})^{-1}\), since
\(\cL_T/\log\log T\to1\) and \(V\le\pi V\). This proves the tail bounds. Likewise
\(\{|S-v|\le h\}=\{|X_\chi-\pi v|\le\pi h\}\), and
Theorem~\ref{thm:local-window} below, at slope \(\pi K+1\) and
half-width \(\pi h\), gives the window bounds after the same passage from \(\cL_T\) to
\(\log\log T\); for the lower bound take
\(\pi h_1\) to be its threshold \(h_0\) and use monotonicity in \(h\).
\end{proof}

For a fixed character and each fixed \(A\), the two-sided order in
\eqref{eq:uniform} was known from Inoue--Li for
\(V\le A(\log\log T)^{2/3}\) without any hypothesis, and for
\(V\le A(\log\log T)^{2/3}(\log\log\log T)^{1/3}\) under RH; see
Section~\ref{sec:related}. On the logarithmic scale, the fixed-slope
asymptotic \(\log\PT\{\eps X_\chi>\kappa\log\log T\}\sim-\kappa^2\log\log T\)
was proved for \(\zeta\) under RH by Wahl~\cite{Wahl}, and it also
follows from the results of Inoue--Li under RH. The zero theorem now
supplies the hypothesis of the results stated under RH. The new range is
\(V\asymp\log\log T\), where the earlier bounds leave an unbounded
factor. The same method gives fixed-width windows
(Theorem~\ref{thm:local-window}) and exponential moments
(Corollary~\ref{cor:moments}).

\subsection{The one-sided inequality}

The analytic input is a \emph{one-sided} bound of
the argument by a single Dirichlet polynomial over primes and prime
squares, for each sign. Let \(C_w\ge1\) and \(C_2,C_0,y_0>0\) be
fixed. For each sign \(\eps\in\{-1,1\}\) let
\(w_\eps\colon[0,\infty)\to\mathbb C\) be a weight with
\begin{equation}\label{eq:weight-axioms}
 w_\eps(0)=-i\eps,\qquad |w_\eps(u)|\le C_w,\qquad
 |w_\eps(u)-w_\eps(0)|\le C_wu,\qquad
 w_\eps(u)=0\quad(u>1).
\end{equation}
For a truncation \(y\ge2\), put
\begin{equation}\label{eq:majorant-polynomial}
 P_{\eps,y}(t)=\Rea\sum_{p\le y}
   \frac{\chi(p)}{p^{1/2+it}}\,
   w_\eps\Bigl(\frac{\log p}{\log y}\Bigr)
 +\Rea\sum_{p^2\le y}
   \frac{\chi(p)^2}{2p^{1+2it}}\,
   w_\eps\Bigl(\frac{2\log p}{\log y}\Bigr).
\end{equation}
The required inequality is
\begin{equation}
 \boxed{\displaystyle
 \eps X_\chi(t)\le P_{\eps,y}(t)
 +C_2\frac{\log(qT)}{\log y}+C_0,
 \qquad \eps\in\{-1,1\},\quad y_0\le y\le T,\quad t\in[T,2T].}
 \tag{OS}\label{eq:majorant}
\end{equation}
The constants are independent of \(T,t,y\). At a given \(T\) we require
\textup{(OS)} for every eligible \(y\), almost everywhere in \(t\).
Each proof uses finitely many truncations at that \(T\), so their
exceptional sets have null union. Applying \textup{(OS)} with \(-\eps\)
also gives a lower bound for \(\eps X_\chi\).

\begin{theorem}[Sharp tails under the one-sided inequality]\label{thm:uniform}
Fix \(\theta_0>0\), \(K>0\), and data \(C_w,C_2,C_0,y_0\). There are
\(C<\infty\), \(c>0\) and \(T_0\), depending only on these, such that
for \(T\ge T_0\), every modulus \(1\le q\le T^{\theta_0}\), every
character \(\chi\bmod q\) for which \textup{(OS)} holds at \(T\) with
this data, \(0\le V\le K\cL_T\), and either \(\eps\in\{-1,1\}\), the
bounds \eqref{eq:uniform} hold.
\end{theorem}

\begin{proof}[Proof of Theorem~\ref{cor:grh}]
Write \(d\mid q\) for the inducing conductor. The zero theorem of
Alp\"oge~\cite{Alpoge}, applied at conductor \(d\), puts every zero
of the primitive inducing function in \(0<\Rea s<1\) on
\(\Rea s=1/2\). Its zero-free half-plane statement also covers the
line \(\Rea s=1\), apart from the principal pole. The functional equation
with the conjugate character then excludes zeros of the completed
primitive function on \(\Rea s=0\). At \(s=0\) its gamma factor cancels
any trivial zero of the primitive \(L\)-function. For \(d=1\) we use
the zeta theorem, with
the poles of the completion at \(0\) and \(1\) accounted for separately
in Appendix~\ref{app:majorant}. The formal statements and independent
checks are recorded in~\cite{Sutherland}.
Corollary~\ref{app:cor:all-characters} therefore
gives \textup{(OS)} with a single choice of data for every \(q,\chi\).
Lemma~\ref{app:lem:missing-euler} accounts for primes dividing \(q\)
that are absent from the inducing modulus. Theorem~\ref{thm:uniform}
now applies to the walk with the primes dividing \(q\) omitted.
\end{proof}

For a fixed modulus the statement can be written with
\(\log\log(qT)\) in place of \(\cL_T\).

\begin{corollary}[Fixed modulus]
\label{thm:main}
Fix \(q\) and \(\chi\), assume \textup{(OS)}, and put \(L_T=\log\log(qT)\).
\begin{enumerate}
\item[\textup{(a)}] For every \(K>0\) there are \(C<\infty\) and \(T_0\)
such that, for \(T\ge T_0\), \(0\le V\le KL_T\), and either
\(\eps\in\{-1,1\}\),
\begin{equation}\label{eq:main-upper}
 \PT\{\eps X_\chi>V\}
 \le C\,\frac{\sqrt{L_T}}{V+\sqrt{L_T}}\,e^{-V^2/L_T}.
\end{equation}
\item[\textup{(b)}] For every \(K>0\) there are \(c>0\) and \(T_0\) such
that, for \(T\ge T_0\), \(0\le V\le KL_T\), and either
\(\eps\in\{-1,1\}\),
\begin{equation}\label{eq:main}
 \PT\{\eps X_\chi>V\}\ge c\,\frac{\sqrt{L_T}}{V+\sqrt{L_T}}\,e^{-V^2/L_T}.
\end{equation}
\end{enumerate}
The constants may depend on \(q,\chi,K,C_w,C_2,C_0,y_0\).
\end{corollary}

The proof from Theorem~\ref{thm:uniform} is given in
Section~\ref{lw:sec:tails}. By Corollary~\ref{app:cor:all-characters}
and the zero theorem, \textup{(OS)} holds for every character, so both
parts hold for every fixed \(\chi\).

\begin{corollary}[Exponential moments]\label{cor:moments}
Assume \textup{(OS)}. For every fixed \(k>0\) and \(\eps\in\{-1,1\}\),
\[
 \frac1T\int_T^{2T}e^{2k\eps X_\chi(t)}\,dt\asymp(\log T)^{k^2},
\]
with constants depending on \(k\) and on the data listed in
Corollary~\textup{\ref{thm:main}}. Uniformly in the setting
of Theorem~\textup{\ref{thm:uniform}},
\[
 \frac1T\int_T^{2T}e^{2k\eps X_\chi(t)}\,dt\asymp e^{k^2\cL_T}
 \asymp\Bigl(\frac{\varphi(q)}q\log T\Bigr)^{k^2},
\]
with constants depending only on \(k\), \(\theta_0\) and the data.
With the common data of Corollary~\ref{app:cor:all-characters}, this
holds for all characters modulo \(q\le T^{\theta_0}\), with constants
depending only on \(k\) and \(\theta_0\).
\end{corollary}

The proof, at the end of Section~\ref{lw:sec:tails}, combines
Theorem~\ref{thm:uniform} with a bound for the far tail, obtained
from \textup{(OS)} at a shorter truncation and a high moment of the
resulting polynomial; no other input from the zeros is used. The
comparison with known moment bounds is in Section~\ref{sec:related}.

\begin{theorem}[Fixed-width windows for the argument]\label{thm:local-window}
Fix \(K,\theta_0>0\) and the data of \textup{(OS)}. For each fixed
half-width \(h>0\), there are constants \(C_h,T_h\) such that for
\(T\ge T_h\), every character in Theorem~\ref{thm:uniform}, and every
real \(v\) with \(|v|\le K\cL_T\),
\[
 \PT\{|X_\chi-v|\le h\}\le C_h\,\cL_T^{-1/2}\exp(-v^2/\cL_T).
\]
There are \(h_0,c,T_0>0\), depending only on \(K,\theta_0\) and the
data, such that, for \(T\ge T_0\) and the same \(q,\chi,v\),
\[
 \PT\{|X_\chi-v|\le h_0\}\ge c\,\cL_T^{-1/2}\exp(-v^2/\cL_T).
\]
Thus the two estimates give matching order for each fixed \(h\ge h_0\).
Corollary~\ref{app:cor:all-characters} and the zero theorem used in
Theorem~\ref{cor:grh} give this conclusion for all characters modulo
\(q\le T^{\theta_0}\). The upper-bound constants may also depend on \(h\).
\end{theorem}

At \(v=0\) the probability of a sufficiently wide fixed window has
order \(\cL_T^{-1/2}\), while the tail at zero has order one.
The proof in Section~\ref{sec:local-window} retains the local factor
by tracking exits on both sides and controlling the two signed
corrections on the same terminal window. The bounds are estimates up
to multiplicative constants for each fixed slope range and width.

\subsection{Structure of the proof}

We compare the polynomial in \textup{(OS)} with a sum over primes,
using independent uniform phases as a model. A positive polynomial
controls the correction on each capped unit window of the prime sum.
The primes in the first block supply the inverse-square-root factor,
and a mean-square estimate for rational frequencies transfers the bound
to the height average. A signed polynomial gives a matching lower
bound for the mass of a capped unit window.

To pass to the argument itself, we track the first time the prime sum
leaves either side of the target interval. For the lower bound we
control both signed corrections on the same terminal window.
Summing the resulting fixed-width estimates gives the tail prefactor
\(\sqrt{\cL_T}/(V+\sqrt{\cL_T})\), including at \(V=0\).
A high moment at a shorter prime cutoff handles the range beyond the
last window. Throughout, the variance clock \(\cL_T\) omits the primes
dividing \(q\); the block sizes and degree budgets use \(\log\log T\).

We use the classical prime estimates
\begin{equation}\label{eq:prime-sums}
 \sum_{p\le z}\frac1p=\log\log z+O(1),\qquad
 \sum_{p\le z}\frac{\log p}{p}\ll\log z,\qquad
 \sum_{p\le z}\frac{(\log p)^2}{p}\ll(\log z)^2.
\end{equation}
The variance clock (Lemma~\ref{lem:clock}), used in both halves, and the
lower local-limit argument use the stronger Mertens formula
with a constant and \(o(1)\) remainder. These estimates are
unconditional. The constants in the intermediate comparisons may
depend on their displayed fixed parameters; uniformity in moment
degrees, scale indices, endpoints and coefficients is asserted only
under the corresponding hypotheses.

\section{Relation to previous results}\label{sec:related}

For the comparisons in this section put $W_0=\log\log T$.
For a fixed modulus $q$, replacing $W_0$ by $W=\log\log(qT)$ changes
$\sqrt{W_0}\,(V+\sqrt{W_0})^{-1}\exp(-V^2/W_0)$ by a factor $1+o(1)$
uniformly for $0\le V\le KW$.
On this whole range the conclusion of
Corollary~\ref{thm:main} is a sharp order estimate:
it includes both the quadratic exponent and the prefactor
$\sqrt W/(V+\sqrt W)$, but does not specify a leading constant or a
relative asymptotic.

Each result below is quoted under the hypothesis stated by its authors.
Where that hypothesis is RH for $\zeta$ or for a Dirichlet
$L$-function, the zero theorem \cite{Alpoge} now supplies it; we keep
the words ``under RH'' to record what each proof uses.

\paragraph{Argument tails.}
Radziwi\l\l's Theorem~1 and Remarks~(1)--(3) in
\cite{Radziwill} give Gaussian relative asymptotics in a moderate
deviation range, including the argument.  In the physical threshold
variable, the stated ranges are $V\ll W_0^{3/5-\varepsilon}$
unconditionally and $V\leq W_0^{2/3-\varepsilon}$ for the argument
under RH.  These ranges precede the scale $V\asymp W_0$ considered here.

Under RH, Wahl \cite[Corollary~1.2]{Wahl} proves that, for $t$
uniform on $[T,2T]$, the family
$(W_0/2)^{-1}\operatorname{Im}\log\zeta(\tfrac12+it)$ satisfies the
large deviation principle with speed $(W_0/2)^{-1}$ and rate function
$h^2/2$; his displayed consequence concerns the upper tail. Applied to
the sets $(2\kappa,\infty)$ and $(-\infty,-2\kappa)$, on whose
interiors and closures the rate function has the same infimum, the
principle gives, for each fixed $\kappa>0$ and either sign,
\[
 \log\Bigl(\frac1T\operatorname{meas}
 \{t\in[T,2T]:\pm\operatorname{Im}\log\zeta(\tfrac12+it)>\kappa W_0\}\Bigr)
 =-(\kappa^2+o(1))\,W_0 .
\]
This fixes the exponent at each fixed positive slope on the logarithmic
scale and leaves a factor $e^{o(W_0)}$ unspecified.

Inoue--Li \cite[Theorem~2.3 and Corollary~2.2]{InoueLi}
provide the closest direct comparisons.  For $F=\zeta$ or a fixed
primitive Dirichlet $L$-function, their variance parameter is $n_F=1$.
Under the corresponding RH, Corollary~2.2, with
$\theta=\pm\pi/2$, gives
\begin{gather*}
 \frac1T\operatorname{meas}
 \{t\in[T,2T]:\pm\operatorname{Im}\log F(\tfrac12+it)>V\}
 \ll_{A,F}\frac{\sqrt{W_0}}{V}e^{-V^2/W_0},\\
 \sqrt{W_0}\leq V\leq A W_0^{2/3}(\log W_0)^{1/3}.
\end{gather*}
On this range the upper bound \eqref{eq:main-upper} has the same
shape as their Corollary~2.2 and is not new; part~(a) of
Corollary~\ref{thm:main} extends that shape to $V\le KW_0$.
Unconditionally, their Corollary~2.1 gives, for
$0\le V\le a_3W_0$ with a small constant $a_3$, upper and lower bounds
for these tails of the form
$(1+V/\sqrt{W_0})^{-1}\exp\{-(V^2/W_0)(1+O(V/W_0))\}$, the upper bound
with an additional term of relative size $O(W_0^{-2})$ in the
prefactor. For $V\le AW_0^{2/3}$, with $A$ fixed, the factor
$\exp\{O(V^3/W_0^2)\}$ is bounded, so in that range the two-sided order
$\sqrt{W_0}\,(V+\sqrt{W_0})^{-1}e^{-V^2/W_0}$ is already known without
any hypothesis; at $V\asymp W_0$ it is $e^{O(W_0)}$.
Their Theorem~2.3 also reaches every fixed positive slope.
After its normalization $z=V/\sqrt{W_0/2}$, the error
$O(z^3/(\sqrt{W_0}\log z))$ is $O_I(W_0/\log W_0)$
when $V/W_0\in I\Subset(0,\infty)$.
Its upper and lower bounds thus allow a growing exponential
factor relative to the sharp order obtained here. For
$3\sqrt{W_0/2}\le V\le AW_0^{2/3}(\log W_0)^{1/3}$, with $A$ fixed, the same error is
$O_A(1)$, and its range conditions hold for large $T$; its lower bound
is stated for $\theta\in[\pi/2,3\pi/2]$, which covers both signs.
Together with Corollary~2.1 for small $V$, Theorem~2.3 therefore gives
the two-sided order above under RH for
$0\le V\le AW_0^{2/3}(\log W_0)^{1/3}$.

For $\zeta$, Inoue \cite[Corollary~1]{Inoue}
proves under RH, for $k\geq0$ and $T$ above the stated
$k$-dependent threshold, the bound
\[
 \frac1T\int_T^{2T}
 \exp\!\bigl(2k\operatorname{Re}(e^{-i\theta}
                    \log\zeta(\tfrac12+it))\bigr)\,dt
 \ll e^{e^{O(k)}}e^{k^2W_0},
 \qquad -\pi/2\leq\theta\leq\pi/2.
\]
At fixed positive slopes, optimizing Markov's inequality with
$k=V/W_0$ gives $O_I(e^{-V^2/W_0})$ for each argument tail.
This deduction alone does not give the factor $W_0^{-1/2}$.

For a primitive non-principal character to a fixed odd squarefree
modulus $q\ge3$, the author \cite[Theorem~1.1]{HughesP1} proved under
GRH that $\frac1T\int_T^{2T}|X_\chi(t)|^{2k}\,dt\le(C_KkW)^k$ for
every integer $1\le k\le KW$, with an unspecified constant $C_K$, and
deduced by Markov's inequality \cite[Corollary~5.1]{HughesP1} that
\[
 \frac1T\operatorname{meas}\{t\in[T,2T]:|X_\chi(t)|>V\}
 \le\exp\Bigl(-\frac{V^2}{2eC_KW}\Bigr),
 \qquad \sqrt{eC_KW}\le V\le\sqrt{eC_KK}\,W .
\]
The remark closing that paper leaves the sharp coefficient $1$ of
$V^2/W$ and the missing prefactor to future work. Corollary~\ref{thm:main}
supplies both for each sign separately, together with the matching
lower bound.

\paragraph{Fixed-width intervals.}
For a fixed primitive Dirichlet \(L\)-function, the interval estimate of Selberg
recorded by Hejhal~\cite[p.~224, (3)]{Hejhal} applies directly to the
imaginary part of the logarithm. Its normalized additive error is
\(O((\log Q)^2/\sqrt Q)\), uniformly in the interval endpoints, with
\(Q=\sum_{p\le T}|\chi(p)|^2/p\) for a primitive Dirichlet character.
At the origin the Gaussian term for a fixed interval is of order
\(Q^{-1/2}\), so the stated error bound exceeds the scale of the term.
For \(\zeta\), Ivi\'c and Korolev~\cite[Theorem~1, (2.1)]{IvicKorolev}
give a smaller endpoint-uniform bound for the additive error on \([T,2T]\), of order
\((\log\log\log T)/\sqrt{\log\log T}\). This still exceeds the order
of a fixed window at the origin by a factor of \(\log\log\log T\).
Kowalski and Nikeghbali~\cite[Theorem~2 and Corollary~9]{KowalskiNikeghbali}
give local estimates for finite Euler products and a criterion for
density of zeta values under a quantitative characteristic-function
conjecture.

There is also a direct comparison with the tail bounds in this paper.
Fix \(0<\delta<K\). For \(a=|v|\in[\delta Q,KQ]\), subtract the signed
tails at \(a-h\) and \(a+h\), using slope range \(K+1\).
The ratio of the upper estimate at the second endpoint to the lower
estimate at the first is \(\ll\exp(-4ah/Q)\), with a constant depending only on
\(K,\theta_0\) and the majorant data.
For sufficiently large fixed \(h\) this yields the lower bound in
Theorem~\ref{thm:local-window} on this smaller range. As \(a/Q\) tends
to zero, the ratio estimate no longer ensures positive mass from
subtraction. Section~\ref{sec:local-window} treats that range directly.

\paragraph{Exponential moments.}
For $\zeta$ under RH, Wahl \cite[Corollary~1.3]{Wahl} proved that,
for every real $h$, $(W_0/2)^{-1}$ times the logarithm of the mean of
$e^{h\operatorname{Im}\log\zeta(\frac12+it)}$ over $t\in[T,2T]$ tends
to $h^2/2$; with $h=\pm2k$, the mean is $(\log T)^{k^2+o(1)}$ for both
signs. Under the same hypothesis, Najnudel \cite{Najnudel} proved that the mean of
$e^{2k\operatorname{Im}\log\zeta(\frac12+it)}$ over $t\in[0,T]$ is
$\ll_{k,\delta}(\log T)^{k^2+\delta}$ for every real $k$ and every
$\delta>0$. With $\theta=\pm\pi/2$, Inoue's bound above removes the
factor $(\log T)^{\delta}$ from the upper bound: it is the upper half
of Corollary~\ref{cor:moments} for $\zeta$, with explicit dependence
on $k$. For each fixed primitive Dirichlet character, Inoue--Li
\cite[Theorem~1.1]{InoueLi} prove under RH, for every fixed $k>0$, every
$\delta>0$ and both signs, that
$\frac1T\int_T^{2T}e^{\pm2k\operatorname{Im}\log L(\frac12+it,\chi)}\,dt$
lies between constant multiples of $(\log T)^{k^2-\delta}$ and
$(\log T)^{k^2+\delta}$. Their Theorem~2.5 \cite{InoueLi}, under RH
with one $L$-function, sharpens this to bounds
$\ll(\log T)^{k^2+Bk^3/\log_3T}$ and $\gg(\log T)^{k^2-Bk^3/\log_3T}$
for both signs, where $\log_3T=\log\log\log T$ and $B$ depends on the
$L$-function; the factor $(\log T)^{Bk^3/\log_3T}$ still tends to
infinity. Corollary~\ref{cor:moments} removes the factors
$(\log T)^{o(1)}$ and $(\log T)^{\pm Bk^3/\log_3T}$; for $\zeta$ only
its lower half goes beyond Inoue's bound. Its proof uses no moment theorem for
$L$-functions from the literature: the range where the signed argument exceeds a large fixed
multiple of $W_0$ is controlled by \textup{(OS)} at a shorter
truncation and a high moment of the resulting prime polynomial, and
Corollary~\ref{thm:main} covers everything below that range. In a limited search we did not find the lower bound of
exact order $(\log T)^{k^2}$ for the argument stated elsewhere. We make
no novelty claim for the corollary beyond this comparison.

Dobner \cite[Theorems~1 and~2]{Dobner} proves lower bounds for
both signs of $S(t)=\pi^{-1}\operatorname{Im}\log\zeta(\tfrac12+it)$,
unconditionally and under RH, over deviation ranges much larger than
$W_0$.  At $V\asymp W_0$, these bounds have Gaussian shape but
leave an unspecified constant in the quadratic exponent.
They therefore do not give the matching sharp order in
Corollary~\ref{thm:main}.

\paragraph{Real-part results and methodological context.}
For $\log|\zeta|$, Arguin--Bailey
\cite[Theorem~1.1]{ABupper} prove the sharp upper bound
unconditionally for $V\sim\alpha W_0$, $0<\alpha<2$,
uniformly on compact slope intervals.
Arguin--Bailey \cite[Theorem~1]{ABlower} prove a sharp lower bound
at $\sigma=\tfrac12+\delta/\log T$, for sufficiently large
$\delta$ depending on the slope.  That theorem is off the critical
line.  Their Proposition~2 compares their particular real-prime
and mollifier barrier event with a random model; it is a precedent
for multiscale comparison, not a theorem for the signed certificate
used here.
Arguin--Creighton \cite[Theorem~1.1]{ArguinCreighton} obtain an
unconditional sharp lower bound for $\log|\zeta|$ on the critical
line, stated on $[0,T]$, for every fixed positive slope.
Arguin--Bailey--Roberts \cite[Theorem~1.1]{ABR2026} prove the
corresponding sharp upper bound on $[T,2T]$ for every fixed positive
slope under RH.  These statements concern the real part of
$\log\zeta$; an argument-tail statement requires its own comparison
with the actual imaginary part.

\paragraph{The conjectured constant.}
For $\log|\zeta|$, Radziwi\l\l{} \cite[Conjecture~2]{Radziwill}
conjectures that at deviations of order $W_0$ the tail probability is
asymptotic to a constant multiple of the Gaussian tail, the constant
being determined by the moment conjecture. The results just cited give
an upper bound of the order of the Gaussian tail on $[T,2T]$,
unconditionally for $V\sim\alpha W_0$ with $0<\alpha<2$ and under RH
for every fixed $\alpha>0$, and a lower bound of the same order on
$[0,T]$, unconditionally for every fixed $\alpha>0$; the constant
remains open. The conjecture is stated for
$\log|\zeta|$ only. The common right side of \eqref{eq:main-upper} and
\eqref{eq:main} is of the order of the Gaussian tail
$\int_{\Delta}^{\infty}e^{-u^2/2}\,du/\sqrt{2\pi}$ with
$\Delta=V/\sqrt{L_T/2}$, uniformly for $0\le V\le KL_T$. For the
argument, Corollary~\ref{thm:main} therefore gives the order of
magnitude of the analogous statement, and nothing about its constant.

\paragraph{The analytic input.}
Selberg-type inequalities bound the argument on both sides by a
Dirichlet polynomial plus the absolute value of a second polynomial.
For $\log|\zeta|$, Soundararajan's method and the recursive scheme of
\cite{ABR2026} use instead a one-sided bound by a single polynomial.
Hypothesis \textup{(OS)} is the analogue for the argument, for each
sign, and for $\zeta$ it is already available in the literature in
two forms. Inoue's Proposition~1 \cite{Inoue} is stated without RH;
under RH and for $\cos\theta\ge0$ it bounds
$\operatorname{Re}(e^{-i\theta}\log\zeta(\tfrac12+it))$ by a
polynomial over primes and prime squares, with an error consisting of
the signed real part of a second prime polynomial and a multiple of
$\log t/\log X$; at $\theta=\pm\pi/2$ this gives \textup{(OS)} for
$\zeta$ with explicit piecewise weights. Replacing Inoue's square sum by
the one in \eqref{eq:majorant-polynomial} costs
$O(\sum_{p\le\sqrt X}\log p/(p\log X)+\sum_{\sqrt X<p\le X}1/p)=O(1)$.
Indeed the two square weights agree at zero and differ by
$O(\log p/\log X)$ on the common range. Our lower bound uses
\textup{(OS)} for the opposite sign. An inequality of this type for
$-\operatorname{Re}(e^{-i\theta}\log L(\tfrac12+it,\chi))$ cannot hold
when $\cos\theta>0$: at a zero on the critical line,
$-\cos\theta\log|L(\tfrac12+it,\chi)|$ tends to $+\infty$, while the
argument stays bounded nearby and the truncated polynomial is
continuous in $t$. For this reason the method gives both bounds only
for the argument. Carneiro--Chandee--Milinovich
\cite{CCM15}, in their proof of the bound for $S(t)$ under RH, pass
through an inequality in which $S(t)$ is bounded by an archimedean
integral minus the prime sum attached to Selberg's majorant of an
interval, together with the pole term $R^+(1/(2i))-1$.
Their formula for $R^+$ and Lemma~\ref{app:lem:H}(b), applied at
$\Delta(t\pm i/2)$, bound this term by $O(e^{\pi\Delta}/(\Delta t)^2)=O(1)$
for $\Delta\ge1$ and $y=e^{2\pi\Delta}\le T\le t$.
Evaluating the interval transform by Vaaler's formula then gives, for
$\zeta$, the weights \eqref{app:eq:weights} with $C_2=\pi/2$.
They then estimate the prime sum trivially, as do Goldston--Gonek
\cite{GoldstonGonek}. Zhao \cite{Zhao} keeps the prime sum produced by
such majorants for Dirichlet $L$-functions to a prime modulus and
averages it over the characters. We therefore do not claim
\textup{(OS)} as new. Appendix~\ref{app:majorant} gives a
self-contained proof for each individual primitive Dirichlet character,
complex ones included: it writes the argument as a sum of a fixed odd
function over the zeros, as in the general theorem of \cite{CCM15},
and majorizes that function by Vaaler's construction for functions of
bounded variation \cite[Corollary~12]{Vaaler}, whose Fourier transform
is explicit. Lemma~\ref{app:lem:missing-euler} transfers the resulting
bound to every induced character with absolute constants.

Theorem~\ref{thm:uniform} gives the sharp order uniformly for
\(q\le T^{\theta_0}\), with the primes dividing \(q\) omitted from the
variance. Its constants apply to every
character of the indicated modulus, including induced characters. Each
probability averages over $t$ for one character. Averages over the characters
modulo $q$ form a different family: in it, Arguin and Creighton
\cite[Theorem~1.1]{ArguinCreightonQ} bound the number of even
primitive characters $\chi$ of large modulus $q\not\equiv2\pmod4$ with $\log|L(\frac12,\chi)|>V$,
divided by $\varphi(q)$, from above for $V\sim\alpha\log\log q$ with
$0<\alpha<1$, and Zhao \cite{Zhao}, for a large odd prime $q$ and under GRH,
bounds the mean and the mean square of
$\arg L(\frac12+it,\chi)+\arg L(\frac12+it,\bar\chi)$ over the
non-principal characters $\chi$ modulo $q$.
For $k=2$, the factor $(\varphi(q)/q)^{4}$ in
Corollary~\ref{cor:moments} has the order of the leading coefficient in
Topacogullari's fourth-moment formula for $|L(\frac12+it,\chi)|$,
which is uniform in the primitive character $\chi$
\cite[Theorem~1.1]{Topacogullari}. That coefficient
equals $\frac{1}{2\pi^2}(\varphi(q)/q)^4\prod_{p\mid q}(1-p^{-2})^{-1}$.

\section{One majorant and one correction polynomial}
\label{sec:foundations}

Fix $\eps\in\{-1,1\}$. Hypothesis \textup{(OS)} is used twice: with the
sign $\eps$ as an upper bound for $\eps X_\chi$, and with the sign
$-\eps$ as a lower bound. Put
\[
 P^{+}_{y}=P_{\eps,y},\qquad P^{-}_{y}=-P_{-\eps,y}.
\]
Both are real parts of polynomials of the shape
\eqref{eq:majorant-polynomial}, with weights $w_\eps$ and $-w_{-\eps}$
respectively, and by \eqref{eq:weight-axioms} both weights take the
value $-i\eps$ at $0$. Fix $\theta_0>0$; throughout
Sections~\ref{sec:foundations}--\ref{lw:sec:tails} the modulus satisfies
$1\le q\le T^{\theta_0}$, and $\chi\bmod q$ may vary with $T$. A fixed
modulus satisfies this for large $T$ with $\theta_0=1$. Set
\[
 B=(1+\theta_0)C_2+C_0 .
\]
For $y=T^{\beta_y}$ with $\beta_y\le1$ and $y\ge\max\{2,y_0\}$ we have
$\log(qT)\le(1+\theta_0)\log T$ and $1\le1/\beta_y$, so
\textup{(OS)} gives
\begin{equation}\label{eq:majorant-final}
 P^{-}_{y}-\frac B{\beta_y}\le \eps X_\chi\le P^{+}_{y}+\frac B{\beta_y}.
\end{equation}
The two outer polynomials are different; no bound for their
difference is used.

For a reference cutoff $z$, define
\begin{equation}\label{eq:reference}
 S_z(t)=\eps\Ima\sum_{\substack{p\le z\\p\nmid q}}
 \frac{\chi(p)}{p^{1/2+it}}
 =\Rea\sum_{\substack{p\le z\\p\nmid q}}
 \frac{-i\eps\,\chi(p)}{p^{1/2+it}}.
\end{equation}
The reference walk contains only first powers of primes, with the
unit weight $-i\eps$.

\begin{lemma}[Coefficient bounds]\label{lem:coefficients}
Suppose $y\ge z\ge2$, and let $W_{y,z}$ denote either
$P^{+}_{y}-S_z$ or $P^{-}_{y}-S_z$.
Write the contribution at base $p$ as
$W_p(t)=\Rea\sum_{k\ge1}a_{p,k}p^{-kit}$ and set
$m_p=\sum_{k\ge1}|a_{p,k}|$.  If $Y_p$ is the base-$p$ summand of
\eqref{eq:reference}, then
\begin{align}
 \sum_p m_p^2&\ll 1+\log\frac{\log y}{\log z},\label{eq:weight-mass}\\
 \sum_{p\le z}\|Y_p\|_\infty m_p&\ll1+\frac{\log z}{\log y}\ll1.
 \label{eq:weight-overlap}
\end{align}
The implied constants depend only on $C_w$. The polynomial $W_{y,z}$
has no constant term, and all its positive and negative frequencies
are supported on integers at most $y$.
\end{lemma}

\begin{proof}
Let $w$ be the weight of the polynomial in question, so $w(0)=-i\eps$,
and write $u_p=\log p/\log y$. Only $k=1,2$ occur. If $p\mid q$ then
$\chi(p)=0$ and $W_p=0$. For $p\le z$ with $p\nmid q$ the reference
walk removes the value $w(0)$ from the first-power coefficient, so by
\eqref{eq:weight-axioms}
\[
 a_{p,1}=\frac{\chi(p)}{\sqrt p}\bigl(w(u_p)-w(0)\bigr),\qquad
 |a_{p,1}|\le C_w\frac{\log p}{\log y}\,p^{-1/2},\qquad
 |a_{p,2}|\le\frac{C_w}{2p}.
\]
For $z<p\le y$ nothing is removed and
$|a_{p,1}|\le C_wp^{-1/2}$, $|a_{p,2}|\le C_w/(2p)$; for $p>y$ the
weight vanishes. Hence
\[
 m_p\le C_w\Bigl(\frac{\log p}{\log y}p^{-1/2}+p^{-1}\Bigr)\quad(p\le z),
 \qquad m_p\le2C_wp^{-1/2}\quad(z<p\le y),
 \qquad m_p=0\quad(p>y).
\]
By \eqref{eq:prime-sums} the low-prime square mass is
$O(1+(\log z/\log y)^2)$, and its overlap with
$\|Y_p\|_\infty=p^{-1/2}$ is
$O(1+\log z/\log y)$, since $\sum_pp^{-3/2}<\infty$.
The high-prime square mass is at most
$4C_w^2\sum_{z<p\le y}1/p\ll1+\log(\log y/\log z)$, which proves
\eqref{eq:weight-mass}. The support assertion is read off from
\eqref{eq:majorant-polynomial}.
\end{proof}

The correction $W_{y,z}$ will remain a single polynomial.  Its
shared-prime and future-prime terms are never separated in the crossing
argument.

\section{One uniform local comparison}
\label{sec:comparison}

\subsection{The reference walk and independent model}

Put $n=\log\log T$, fix $\gamma=1/2$, and let
\begin{equation}\label{eq:ladder}
 \beta_i=\frac{e^{i-1}}{\sqrt{\log T}},\quad
 T_i=T^{\beta_i},\quad t_i=\log\log T_i,\quad
 x_i=\log\beta_i^{-1},\quad c_i=\beta_i^{-\gamma}.
\end{equation}
Write $\mathcal P_i$ for the primes in block $i$. The first block
contains all primes up to $T_1$, and later blocks have $T_{i-1}<p\le T_i$. Set
\[
 Y_i=S_{T_i}-S_{T_{i-1}},\quad S_\ell=\sum_{i\le\ell}Y_i,
 \qquad \nu_1=t_1+1,\quad\nu_i=1\ (i\ge2),
\]
where $S_{T_0}=0$.  In the model, replace $p^{-it}$ by independent random
variables $U_p$, uniform on the unit circle.  Every power of the same
prime uses the corresponding power of the same $U_p$.
We write a superscript $*$ when needed to
distinguish this model from the arithmetic polynomial.

The primes dividing $q$ are absent from the walk. Its variance is
therefore measured by the \emph{variance clock}
\begin{equation}\label{eq:clock}
 \tilde t_\ell=\sum_{\substack{p\le T_\ell\\p\nmid q}}\frac1p,
 \qquad
 \cL_T=\sum_{\substack{p\le T\\p\nmid q}}\frac1p ,
\end{equation}
while the geometric clock $t_\ell$ continues to index the ladder,
degrees, caps and heights. Mertens' formula with a constant
\cite{Mertens},
\begin{equation}\label{lw:eq:mertens}
 \sum_{p\le x}p^{-1}=\log\log x+B_{\rm M}+o(1),
\end{equation}
is unconditional.

\begin{lemma}[The variance clock]\label{lem:clock}
Uniformly for $1\le q\le T^{\theta_0}$,
\begin{enumerate}
\item[\textup{(i)}] $\sum_{p\mid q}p^{-1}\le\log n+O(1)$;
\item[\textup{(ii)}] $\sum_{p\mid q,\ p>T_1}p^{-1}
\le\theta_0\sqrt{\log T}\,e^{-\sqrt{\log T}}$;
\item[\textup{(iii)}] for every $\ell$ with $\beta_\ell\le1$,
$|\tilde t_\ell-t_\ell|\le\log n+O(1)$ and
$\cL_T-\tilde t_\ell=x_\ell+o(1)$; for $\ell\ge2$ also
$\tilde t_\ell-\tilde t_{\ell-1}=1+o(1)$.
\end{enumerate}
The $o(1)$ and $O(1)$ depend only on $\theta_0$. In particular
$\tilde t_\ell/t_\ell\to1$ uniformly in $\ell$, since $t_\ell\ge t_1=n/2$.
\end{lemma}

\begin{proof}
(i) The primes $p\mid q$ with $p\le\log T$ contribute at most
$\sum_{p\le\log T}p^{-1}=\log n+O(1)$. There are at most
$\log q/\log\log T\le\theta_0\log T/\log\log T$ prime factors
$p>\log T$, each contributing less than $1/\log T$; together less
than $\theta_0/n$.
(ii) At most $\log q/\log T_1\le\theta_0\sqrt{\log T}$ prime factors
exceed $T_1=e^{\sqrt{\log T}}$, each contributing less than $1/T_1$.
(iii) By \eqref{lw:eq:mertens}, $\sum_{p\le T_\ell}p^{-1}-t_\ell=O(1)$,
and removing the primes dividing $q$ costs at most (i). For the other
two assertions, the primes in $(T_\ell,T]$, respectively in
$(T_{\ell-1},T_\ell]$, have reciprocal sum
$\log\log T-\log\log T_\ell+o(1)=x_\ell+o(1)$, respectively
$t_\ell-t_{\ell-1}+o(1)=1+o(1)$, by \eqref{lw:eq:mertens} and
$T_{\ell-1}\ge T_1\to\infty$; by (ii) the primes dividing $q$ change
these sums by $o(1)$.
\end{proof}

Each prime contribution has mean zero and
\[
 \|Y_p\|_\infty=p^{-1/2},\qquad
 \Var(Y_p)=\frac1{2p}\quad(p\nmid q).
\]
In particular
\begin{equation}\label{eq:block-mass}
 \sum_{p\text{ in block }i}\|Y_p\|_\infty^2\ll\nu_i,
 \qquad
 \sum_{p\le T_\ell}\Var(Y_p)=\frac12\tilde t_\ell.
\end{equation}
For every bounded real $h$,
\begin{equation}\label{eq:reference-mgf}
 \log\E e^{hS_\ell^*}=\frac{h^2\tilde t_\ell}{4}+O_h(1),
\end{equation}
uniformly on compact $h$-ranges and in $q$ and $\chi$.  Indeed,
symmetry and Taylor expansion give
$\log\E e^{hY_p}=h^2/(4p)+O_R(p^{-2})$ uniformly for $|h|\le R$; sum
this identity over the primes $p\le T_\ell$ with $p\nmid q$.  The
remainder is one global $O(1)$, independent of the number of blocks
and of $q$. Below, $\tilde t_\ell$ appears exactly where the clock
multiplies a tilt or a slope (exponents, barriers, thresholds); the
inverse square roots, ladder, degrees and heights keep $t_\ell$, which
is legitimate by Lemma~\ref{lem:clock}(iii).

We will use the following elementary consequences of boundedness. If
independent real variables $Z_j$ are supported in intervals of length
at most $2b_j$, then
\begin{equation}\label{eq:hoeffding}
 \E e^{s\sum_j(Z_j-\E Z_j)}
 \le e^{s^2\sum_j b_j^2/2},\qquad
 \Big\|\sum_j(Z_j-\E Z_j)\Big\|_{2k}
 \ll\sqrt{k\sum_j b_j^2}.
\end{equation}
For completeness, the second derivative of $\log\E e^{sZ_j}$ is
the variance under the $s$-tilt. This variance is at most $b_j^2$:
subtract the midpoint of the support interval and bound its second
moment. Integrating twice and using independence proves the exponential
bound. Its Chernoff estimate, integrated against $2kv^{2k-1}\,dv$,
gives the moment bound. This is the bounded-variable inequality
associated with Hoeffding \cite{Hoeffding}. Both statements remain
valid after any product tilt, with centering under that tilt, because
the support intervals do not change.

\subsection{Rational-frequency transfer}

A finite Laurent polynomial has height $(A,B)$ if its reduced rational
frequencies $a/b$ satisfy $a\le A$ and $b\le B$.  Its symmetric height is
$H$ when $A,B\le H$.

\begin{lemma}[Mean-square transfer]\label{lem:mv}
If $F(t)=\sum_\rho d_\rho\rho^{-it}$ has height $(A,B)$, then
\[
 \left|\ET|F|^2-\E|F^*|^2\right|
 \le C\frac{AB}{T}\E|F^*|^2.
\]
The constant is independent of all coefficients and degrees.
\end{lemma}

\begin{proof}
For distinct reduced rationals $a/b,c/d$ in this range,
\[
 |\log(a/b)-\log(c/d)|
 =|\log(ad)-\log(bc)|\ge\frac1{AB}.
\]
Expand the integral of $|F|^2$. Each off-diagonal term is the difference
of its values at $T$ and $2T$, divided by its frequency difference.
The Montgomery--Vaughan separated-frequency Hilbert inequality
\cite[Theorem~2, (1.6), and Corollary~2, (1.8)]{MV1974},
applied to the two endpoint forms with coefficients
$d_\rho e^{-iT\log\rho}$ and $d_\rho e^{-2iT\log\rho}$, bounds this
off-diagonal integral by $CAB\sum_\rho|d_\rho|^2$. Unique
factorization and orthogonality of distinct rational monomials give
$\E|F^*|^2=\sum_\rho|d_\rho|^2$.  Division by $T$ proves the result.
\end{proof}

Thus a symmetric height $H$ costs $H^2$, not $H$, in mean square.

\subsection{Statement of the comparison}

Let $W$ be any finite real Laurent polynomial of the form
\[
 W^*=\sum_p W_p,
 \qquad W_p=W_p(U_p),\qquad \E W_p=0,\qquad |W_p|\le m_p,
\]
where each $W_p$ uses only powers of its base prime.  There is no bound
on its harmonic degree; for example one may take
$W_p=\Rea\sum_{k\ge1}a_{p,k}U_p^k$ and $m_p=\sum_k|a_{p,k}|$.
Primes both below and above $T_\ell$ are allowed.  Take any proxy
$\Sigma\ge\sum_pm_p^2$ and suppose
\begin{equation}\label{eq:absolute-overlap}
 \sum_{p\le T_\ell}\|Y_p\|_\infty m_p\le M\sqrt\Sigma
\end{equation}
for a fixed $M$.  This absolute overlap condition is stronger than a
covariance bound but immediately fits all the weights used here.  No
individual coefficient smallness hypothesis is imposed on $W$.

For a fixed $K_0\ge1$, put
\[
 B_1=K_0(t_1+c_1),\qquad B_i=K_0c_i\quad(i\ge2).
\]
Let $\cC$ be any measurable event on which
\begin{equation}\label{eq:caps}
 |Y_i|\le B_i\qquad(1\le i\le\ell).
\end{equation}
Fix a compact range $|\lambda|\le R$ and choose
\begin{equation}\label{eq:degrees}
 N_i=\left\lceil A_{\rm cert}(\nu_i+c_i)\right\rceil,
\end{equation}
where $A_{\rm cert}\ge1$ is chosen once from the explicit sufficient
Taylor calibrations below. It depends only on $R,K_0,M$, before any
terminal cutoff is chosen. No later argument enlarges this constant.

\begin{proposition}[Uniform unit-window comparison]\label{prop:comparison}
Suppose $\beta_\ell\le\beta_*<1$ and $r$ is an integer with
$1\le r\le c_\ell$.  Suppose $W$ has symmetric Laurent height at most
$T^{h_W}$, and, for a fixed $\eta>0$,
\begin{equation}\label{eq:height-budget}
 2\sum_{i\le\ell}\beta_iN_i+2rh_W\le1-\eta.
\end{equation}
Then there is $T_0$, depending only on the fixed parameters listed at
the end of this statement, such that for $T\ge T_0$, uniformly for
every real $u$,
\begin{equation}\label{eq:comparison}
 \ET\left[|W|^{2r}\one_{\cC}\one_{\{S_\ell\in[u,u+1]\}}\right]
 \ll \frac{e^{\lambda^2\tilde t_\ell/4-\lambda u}}{\sqrt{t_\ell}}
       (Cr\Sigma)^r.
\end{equation}
There is also an unweighted version, with the final moment factor
replaced by $1$ and the term $2rh_W$ omitted from
\eqref{eq:height-budget}.
The constants depend only on the fixed parameters
$R,K_0,M,\theta_0,\eta$; they are independent of $u,r,\ell$, of $q$
and $\chi$, of the harmonic degrees, and of the coefficients of $W$
subject to the stated bounds.
\end{proposition}

We prove this proposition in full.  There is no local central limit
theorem or fine-cell estimate in its hypotheses.

\subsection{Weighted Taylor approximation}

Write $E_N(z)=\sum_{j=0}^Nz^j/j!$.  For a block $i$, write $W_i$ for the
part of $W$ on that block, and put
\[
 Z_i(h,s)=\E e^{hY_i+sW_i}.
\]
For this block, $\Pp_{h,s}$ denotes the probability measure with
density $Z_i(h,s)^{-1}e^{hY_i+sW_i}$ relative to the independent model;
$\E_{h,s}$ denotes its expectation. These tilts factor over primes.
If $\Sigma>0$, take $s=\pm\sqrt{r/\Sigma}$.  The case $\Sigma=0$ in
the weighted proposition is trivial.

\begin{lemma}[Weighted Taylor error]\label{lem:taylor}
For $|\xi|\le1$, $|\lambda|\le R$, and the degrees
\eqref{eq:degrees}, uniformly in $i\le\ell$,
\begin{equation}\label{eq:taylor-error}
 \left\|E_{N_i}((\lambda/2+i\xi)Y_i)
       -e^{(\lambda/2+i\xi)Y_i}\right\|_{L^2(e^{sW_i}\,d\Pp)}
 \le e^{-N_i}Z_i(\lambda,s)^{1/2}.
\end{equation}
The same estimate holds when $r=0$ and $s=0$.
\end{lemma}

\begin{proof}
Every tilt in question factors over primes.  Differentiate the mean of
$Y_i$ along a tilt by $sW_i+hY_i$.  The elementary bound
$|\Cov(U,V)|\le4\|U\|_\infty\|V\|_\infty$ gives, for bounded $h$,
\begin{equation}\label{eq:tilted-block-mean}
 |\E_{h,s}Y_i|
 \ll\nu_i+|s|\sum_{p\text{ in block }i}\|Y_p\|_\infty m_p
 \ll\nu_i+\sqrt r.
\end{equation}
Here the $h$-tilt alone changes the mean by $O(\nu_i)$, using
\eqref{eq:block-mass}; the second bound uses
\eqref{eq:absolute-overlap}.  Independence and
\eqref{eq:hoeffding} under the product tilt now imply
\begin{equation}\label{eq:tilted-block-moments}
 \|Y_i\|_{L^{2k}(\Pp_{h,s})}
 \ll\nu_i+\sqrt r+\sqrt{k\nu_i}.
\end{equation}
Integration of \eqref{eq:tilted-block-mean} in $h$ also gives
\begin{equation}\label{eq:normalizer-ratio}
 \frac{Z_i(0,s)}{Z_i(\lambda,s)}
 \le\exp(C(\nu_i+\sqrt r)).
\end{equation}

Let $k=N_i+1$ and put $\rho_0=\sqrt{R^2/4+1}$. The integral Taylor
remainder satisfies
\[
 |E_{N_i}((\lambda/2+i\xi)Y_i)-e^{(\lambda/2+i\xi)Y_i}|
 \le\frac{(\rho_0|Y_i|)^k}{k!}(1+e^{\lambda Y_i/2}).
\]
Minkowski's inequality, \eqref{eq:tilted-block-moments}, and
\eqref{eq:normalizer-ratio} bound its weighted $L^2$ norm, divided by
$Z_i(\lambda,s)^{1/2}$, by
\begin{equation}\label{eq:weighted-taylor-majorant}
 \exp(C(\nu_i+\sqrt r))
 \left[\frac{C(\nu_i+\sqrt r+\sqrt{k\nu_i})}{k}\right]^k.
\end{equation}
Here and in the positive majorants below, fix one constant $C\ge1$
depending only on $R,M$ that makes these estimates valid. This is a
proof constant, distinct from the weight constant $C_w$.
We used $k!\ge(k/e)^k$ and absorbed $\rho_0$ into $C$.

Write $q_i=\nu_i+c_i$. Since $\sqrt r\le c_i$ and
$k\ge A_{\rm cert}q_i$, the expression in
\eqref{eq:weighted-taylor-majorant} is at most
\[
 \left\{
 e^{C/A_{\rm cert}}C
       (A_{\rm cert}^{-1}+A_{\rm cert}^{-1/2})
 \right\}^{k}.
\]
Choose $A_{\rm cert}$ so that the expression in braces is at most
$e^{-1}$. This proves \eqref{eq:taylor-error}, uniformly in $r$.
The same choice also gives the following explicit unweighted
positive-majorant calibration, with the very same $C$:
\begin{equation}\label{eq:calibrated-moment}
 \boxed{\quad
 e^{C\nu_i}
 \left[C\frac{\nu_i+\sqrt{(N_i+1)\nu_i}}{N_i+1}\right]^{N_i+1}
 \le e^{-N_i}.\quad}
\end{equation}
More generally, for every integer $j\ge1$ and $\rho>0$, the
unweighted integral Taylor majorant obeys
\begin{align}
 &\frac{\rho^j}{j!}
 \left\{Z_i(\lambda,0)^{-1/2}(\E|Y_i|^{2j})^{1/2}
                   +(\E_{\lambda,0}|Y_i|^{2j})^{1/2}\right\}
 \nonumber\\[-2pt]
 &\hspace{25pt}\le
 e^{C\nu_i}\left[
 C\frac{\rho}{\rho_0}\frac{\nu_i+\sqrt{j\nu_i}}j
 \right]^j.
 \label{eq:unweighted-taylor-majorant}
\end{align}
This records a positive sufficient majorant, rather than merely the
smallness of an actual remainder; later degree changes may therefore
use its displayed dependence on $j$ and $\rho$.

We simultaneously calibrate the exact factorial cap bound. Since
$B_i\le K_0q_i$, it is enough to require
\[
 \frac{2e\rho_0K_0}{A_{\rm cert}}
       e^{RK_0/(2A_{\rm cert})}\le e^{-1}.
\]
Indeed $k!\ge(k/e)^k$ and
$1+e^{RB_i/2}\le2e^{RK_0k/(2A_{\rm cert})}$ then give
\begin{equation}\label{eq:calibrated-cap}
 \boxed{\quad
 (1+e^{RB_i/2})
       \frac{(\rho_0B_i)^{N_i+1}}{(N_i+1)!}
 \le e^{-N_i}.\quad}
\end{equation}
Finally, as part of this same initial choice, put
$a_0=\sin(1/2)/(1/2)$ and require
\begin{equation}\label{eq:calibrated-product}
 \sigma_A:=\frac{e^{-A_{\rm cert}}}{1-e^{-A_{\rm cert}/2}},
 \qquad e^{\sigma_A}-1\le a_0/2.
\end{equation}
All three scalar requirements hold once $A_{\rm cert}$ is sufficiently
large. They depend only on the already fixed parameters and hold
before the choice of $\beta_*$. This completes the calibration and
the proof of the lemma.
\end{proof}

\subsection{A pointwise positive certificate}

For $v=u+1/2$, define the finite Laurent polynomial
\begin{equation}\label{eq:certificate}
 F_v=\frac12\int_{-1}^1e^{-i\xi v}
       \prod_{i\le\ell}E_{N_i}((\lambda/2+i\xi)Y_i)\,d\xi.
\end{equation}
Integration changes coefficients but not frequency support.  With
untruncated exponentials this expression is
\[
 e^{\lambda S_\ell/2}\frac{\sin(S_\ell-v)}{S_\ell-v}.
\]
On $S_\ell\in[u,u+1]$, the sine quotient is at least
$a_0=\sin(1/2)/(1/2)>0$. On the caps \eqref{eq:caps}, writing
$z=(\lambda/2+i\xi)Y_i$ and $\rho_0=\sqrt{R^2/4+1}$, the factorial remainder gives
\[
 |E_{N_i}(z)e^{-z}-1|
 \le\frac{(\rho_0|Y_i|)^{N_i+1}}{(N_i+1)!}
             (1+e^{R|Y_i|/2})
 \le e^{-N_i}
\]
by \eqref{eq:calibrated-cap}. This is a relative estimate for either
sign of $Y_i$, uniform for $|\xi|\le1$.
Moreover,
\[
 \sum_{i\le\ell}e^{-N_i}
 \le\sum_{j\ge0}e^{-A_{\rm cert}c_\ell e^{j/2}}
 \le\frac{e^{-A_{\rm cert}}}{1-e^{-A_{\rm cert}/2}}
 =\sigma_A,
\]
where $c_\ell\ge1$ and $e^{j/2}\ge1+j/2$. The relative difference
of the product from its untruncated version is at most
$e^{\sigma_A}-1\le a_0/2$. After integration the absolute error
is at most $(a_0/2)e^{\lambda S_\ell/2}$, since the untruncated
integrand has modulus $e^{\lambda S_\ell/2}$. Hence on the capped window
$|F_v|\ge(a_0/2)e^{\lambda S_\ell/2}$. Since
$e^{\lambda(S_\ell-u)}\ge e^{-R}$ there, this proves
\begin{equation}\label{eq:pointwise-certificate}
 \one_{\cC}\one_{\{S_\ell\in[u,u+1]\}}
 \le C e^{-\lambda u}|F_v|^2.
\end{equation}
No assertion about the approximation outside the caps is required.

The symmetric height of $F_v$ is at most
$T^{\sum_{i\le\ell}\beta_iN_i}$, since $Y_i$ has symmetric height
at most $T_i$.  Consequently \eqref{eq:height-budget} and
Lemma~\ref{lem:mv} give the multiplicative comparison
\begin{equation}\label{eq:certificate-transfer}
 \ET|W^rF_v|^2=(1+O(T^{-\eta}))\E|(W^*)^rF_v^*|^2.
\end{equation}
Shared primes between $W$ and $F_v$ cause no difficulty in this step.

\subsection{The weighted model estimate}

We claim that, for both signs of $s=\sqrt{r/\Sigma}$,
\begin{equation}\label{eq:weighted-model-target}
 \E[e^{sW^*}|F_v^*|^2]
 \ll\frac{e^{\lambda^2\tilde t_\ell/4+Cr}}{\sqrt{t_\ell}}.
\end{equation}
First tilt by $\lambda S_\ell$.  As in
\eqref{eq:tilted-block-mean}, centering under the original measure and
\eqref{eq:absolute-overlap} imply
\[
 |\E_\lambda W^*|\ll\sqrt\Sigma.
\]
Hoeffding's bound under this product tilt and
\eqref{eq:reference-mgf} imply
\begin{equation}\label{eq:joint-mgf}
 Z:=\E e^{\lambda S_\ell^*+sW^*}
 \le\exp\left(\frac{\lambda^2\tilde t_\ell}{4}
       +O(1)+C|s|\sqrt\Sigma+Cs^2\Sigma\right)
 \ll_{R}e^{\lambda^2\tilde t_\ell/4+C_{R,M}r}.
\end{equation}

The first block supplies the density factor.  Call a prime bad if
$|s|m_p>1$.  There are at most $r$ bad primes, because
$s^2\sum m_p^2\le r$.  Discard the primes at most
$P=K_1(r+1)$, with a fixed $K_1\ge\max(64,R^2)$, all
remaining bad primes, and the primes dividing $q$.  By
Lemma~\ref{lem:clock}(i), their reciprocal mass is at most
\[
 \log\log(r+3)+\log n+O(1)+\frac r{K_1(r+1)}.
\]
Since $r\le c_\ell\le c_1=e^{\gamma t_1}$ and $n=2t_1$, the retained
primes in the first block have reciprocal mass at least
$t_1-O(\log(t_1+2))$, hence at least $t_1/2$ for large $T$.

On a retained prime, the one-prime density of the joint tilt
$\lambda Y_p+sW_p$ is bounded above and below by fixed positive
constants.  It follows that
\[
 \Var_{\lambda,s}(Y_p)\gg \Var(Y_p)\gg p^{-1}.
\]
For $|\tau|\le2$, square the modulus of its characteristic function
and introduce an independent copy $Y_p'$.  The bound
$1-\cos x\gg x^2$ on the relevant bounded interval gives
\[
 1-|\E_{\lambda,s}e^{i\tau Y_p}|^2
 =\E_{\lambda,s}[1-\cos(\tau(Y_p-Y_p'))]
 \gg\frac{\tau^2}{p}.
\]
Multiplication over the retained primes, with modulus at most one for
the discarded factors, proves
\begin{equation}\label{eq:first-block-damping}
 |\E_{\lambda,s}e^{i\tau Y_1}|
 \le e^{-ct_1\tau^2}\qquad(|\tau|\le2).
\end{equation}
This deletion argument is why no maximum-coefficient assumption on $W$
is needed.

Expand the square of \eqref{eq:certificate} as a double integral in
$\xi,\zeta\in[-1,1]$.  Independence over prime blocks gives block
kernels
\[
 K_i(\xi,\zeta)=\E\left[e^{sW_i}
 E_{N_i}((\lambda/2+i\xi)Y_i)
 \overline{E_{N_i}((\lambda/2+i\zeta)Y_i)}\right].
\]
By Lemma~\ref{lem:taylor} and \eqref{eq:first-block-damping},
\[
 |K_1(\xi,\zeta)|\le
 Z_1(\lambda,s)\{e^{-ct_1(\xi-\zeta)^2}+O(e^{-N_1})\}.
\]
For each later block, Cauchy--Schwarz and the same Taylor bound give
\[
 |K_i(\xi,\zeta)|\le(1+O(e^{-N_i}))Z_i(\lambda,s).
\]
The future-prime part of $W$ contributes its independent exponential
moment.  The product of all later-block error factors is bounded.
Finally,
\[
 \int_{-1}^1\int_{-1}^1e^{-ct_1(\xi-\zeta)^2}\,d\xi\,d\zeta
 \ll t_1^{-1/2},\qquad e^{-N_1}\ll t_1^{-1/2}.
\]
As $t_1\le t_\ell\le2t_1$, \eqref{eq:joint-mgf} proves
\eqref{eq:weighted-model-target}.  Later-block Taylor errors have been
kept multiplicative; treating them as global additive errors would lose
this density factor.

The elementary inequality
\[
 |w|^{2r}\le\left(\frac{2r}{e|s|}\right)^{2r}
              (e^{sw}+e^{-sw})
\]
and \eqref{eq:weighted-model-target} yield
\[
 \E[|W^*|^{2r}|F_v^*|^2]
 \ll (Cr\Sigma)^r\frac{e^{\lambda^2\tilde t_\ell/4}}{\sqrt{t_\ell}}.
\]
Together with \eqref{eq:pointwise-certificate} and
\eqref{eq:certificate-transfer}, this proves
Proposition~\ref{prop:comparison}.  The unweighted case is its identical
$s=0$ argument without moment extraction.

\section{Parameters and the first block}\label{sec:parameters}

Fix a slope bound $K\ge1$ independently of $T$.
Fix a real tilt range $[-R,R]$ with $R\ge2K$, the cap constant
$K_0=\max(K,4)$, and the coefficient and overlap constants in
Lemma~\ref{lem:coefficients}. The range $[-R,R]$ may be chosen larger
at this initial stage to contain twice any additional fixed compact
slope interval needed for the lower bound. Choose the remaining
constants in the following order:
\begin{enumerate}
\item Choose $A\ge\max(8B,2)$, where $B$ is from
\eqref{eq:majorant-final}.
\item Choose $D\ge100A$ and put $r_\ell=\lfloor c_\ell/D\rfloor$.
\item Fix $A_{\rm cert}$ for Proposition~\ref{prop:comparison}, using
the chosen range $|\lambda|\le R$ and the overlap bounds below.
Require all three explicit sufficient calibrations
\eqref{eq:calibrated-moment}--\eqref{eq:calibrated-product}, as proved
in Lemma~\ref{lem:taylor}.
\item Define the fixed constants
\begin{equation}\label{eq:geometric-constant}
 C_{\rm geo}=\frac2{1-e^{-1}}+\frac1{1-e^{-1/2}},\qquad
 H_0=\max(2D,4K,2),
\end{equation}
and choose a fixed terminal cutoff satisfying the explicit contract
\begin{equation}\label{eq:upper-cutoff}
 0<\beta_*\le
 \min\left\{1,H_0^{-2},
 A^{-2},
                 (32C_{\rm geo}A_{\rm cert})^{-2}\right\}.
\end{equation}
Set
$L=\max\{\ell:\beta_\ell\le\beta_*\}$.
\item Finally increase $T_0$ so that every preceding uniform estimate
holds and $L\ge2$.
\end{enumerate}
In summary, $R$ depends on the slope ranges; $K_0$ on $K$; $A$ and
$D$ on $B$; $A_{\rm cert}$ on $R,K_0,M$; $\beta_*$ on
$H_0,A,A_{\rm cert}$. The lower bound later fixes $a_{\rm LW}$ from
these, and then a threshold $H$ from $a_{\rm LW}$; neither feeds back
into an earlier choice.
This fixes every parameter before $T$ varies. In particular
$c_L\ge H_0$ and $A\beta_*^\gamma\le1$.
The latter ensures that
\begin{equation}\label{eq:adaptive-y}
 y_\ell=T^{A/c_\ell},\qquad \beta_{y,\ell}=A/c_\ell,
\end{equation}
is admissible in \textup{(OS)}, that is, $y_\ell\le T$; it also satisfies
$y_\ell\ge T_\ell$.
Indeed $A\sqrt{\beta_\ell}\ge\beta_\ell$ when $A\ge1$ and
$\beta_\ell\le1$.  The lower cutoff $y_\ell\ge y_0$ is harmless since
$\log y_\ell\ge\log y_1=A(\log T)^{3/4}\to\infty$.

The ceiling in \eqref{eq:degrees} and $A_{\rm cert}\ge1$ give the
following exact geometric estimate:
\begin{equation}\label{eq:certificate-cost}
 \sum_{i\le\ell}\beta_iN_i
 \le A_{\rm cert}\beta_1t_1
                  +C_{\rm geo}A_{\rm cert}\sqrt{\beta_\ell}.
\end{equation}
Indeed the left side is at most
\[
 A_{\rm cert}\beta_1t_1
 +(A_{\rm cert}+1)\sum_{i\le\ell}\beta_i
 +A_{\rm cert}\sum_{i\le\ell}\sqrt{\beta_i},
\]
and the two geometric sums are bounded by
$\beta_\ell/(1-e^{-1})$ and
$\sqrt{\beta_\ell}/(1-e^{-1/2})$, respectively. Use
$\beta_\ell\le1$ to obtain \eqref{eq:certificate-cost}.
The contract \eqref{eq:upper-cutoff} makes twice its right side
less than $1/8$ for sufficiently large $T$.
For either correction $W_\ell=P_{y_\ell}^{\pm}-S_\ell$,
the additional mean-square height is
\begin{equation}\label{eq:weighted-height}
 2r_\ell h_W\le2r_\ell\beta_{y,\ell}
 \le\frac{2A}{D}<\frac18.
\end{equation}
For a reference increment $Y_\ell$ it is at most
$2r_\ell\beta_\ell\le2\sqrt{\beta_\ell}/D$.
Thus every comparison below satisfies
\eqref{eq:height-budget}, for example with $\eta=1/2$.

By Lemma~\ref{lem:coefficients}, the merged weights have proxies
\begin{equation}\label{eq:merged-proxy}
 \Sigma_\ell=C_3 M(x_\ell),\qquad
 M(x)=1+\log A+(1-\gamma)x,
\end{equation}
with $C_3$ fixed sufficiently large.  Their absolute overlap with the
reference walk is $O(1)\ll\sqrt{\Sigma_\ell}$.  A future increment
$Y_\ell$, viewed from level $\ell-1$, has zero overlap and proxy $O(1)$.
All weights therefore fit the same comparison proposition.

Put $G_1=\{|S_1|\le c_1/2\}$.
The entire first bad event can be discarded by an ordinary moment
bound.  Hoeffding's estimate and Lemma~\ref{lem:mv} give
\begin{equation}\label{eq:initial-moment}
 \ET|S_1|^{2r}\le(Crn)^r
 \quad\text{if }2r\beta_1\le\tfrac12.
\end{equation}
The condition holds for $r=r_1$.  Since
$c_1=e^{\gamma n/2}\gg n$, for sufficiently large $T$,
\begin{equation}\label{eq:initial-error}
 \PT(G_1^c)
 =\PT\{|S_1|>c_1/2\}
 \le\left(\frac{Cr_1n}{c_1^2}\right)^{r_1}
 \le\frac{e^{-K^2(n+1)}}{Kn}e^{-c_1}.
\end{equation}
The last inequality holds because
$r_1\asymp c_1/D$ and $\log(c_1/n)\sim\gamma n/2$.

\section{A lower bound for the terminal unit window}\label{lw:sec:statement}

Use the arithmetic walk, independent model, and parameters fixed in
Section~\ref{sec:parameters}. In this part, $\tau\in[T,2T]$ denotes
the arithmetic variable. Thus $n=\log\log T$,
$\beta_i=e^{i-1}/\sqrt{\log T}$, $T_i=T^{\beta_i}$,
$t_i=\log\log T_i$, and $c_i=\beta_i^{-1/2}$. At the terminal index
$L=\max\{i:\beta_i\le\beta_*\}$, write
\[
 t=t_L,\qquad \tilde t=\tilde t_L,\qquad x_L=\log\beta_L^{-1},\qquad x_*=\log\beta_*^{-1},
\]
where $\tilde t_L$ is the variance clock \eqref{eq:clock}.
For all sufficiently large $T$,
\begin{equation}\label{lw:eq:ladder}
 \begin{gathered}
 t_1=n/2,\quad t_i=t_1+i-1,\quad t=n-x_L,\quad
 x_L\in[x_*,x_*+1),\\
 t_1\le t\le2t_1,\qquad c_i=c_Le^{(L-i)/2},\qquad c_1=e^{t_1/2}.
 \end{gathered}
\end{equation}
The same rectangular caps are
\begin{equation}\label{lw:eq:caps}
 B_1=K_0(t_1+c_1),\qquad B_i=K_0c_i\ (i\ge2),\qquad
 \cC_L=\bigcap_{i\le L}\{|Y_i|\le B_i\}.
\end{equation}
Fix $K>0$ and put $J=[0,K]$ before choosing the parameters. The
common real tilt range $[-R,R]$ contains $2J$, so $R\ge2K$. The left
endpoint $0$ is allowed: the proof below uses the slope only through
$|\lambda|\le R$, $\kappa\le K$, and the lower bound
$\log Z\ge\lambda^2\tilde t/4-C_R$, all valid at $\lambda=0$. The degree parameter
$\Ac$, the caps, and the cutoff have already been fixed, with
\[
 \nu_1=t_1+1,\quad\nu_i=1\ (i\ge2),\qquad
 N_i=\lceil\Ac(\nu_i+c_i)\rceil,\qquad
 \Cg=\frac2{1-e^{-1}}+\frac1{1-e^{-1/2}}.
\]
The part of their common cutoff contract used here is
\begin{equation}\label{lw:eq:cutoff}
 \Ac\ge1,\quad K_0\ge4,\qquad
 0<\beta_*\le\min\{400^{-2},(32\Cg\Ac)^{-2}\}.
\end{equation}
The proof below works for every such $K_0$ with the calibrated
$\Ac$. In particular, it does not require a further increase of the
cap constant to $K$.

\begin{theorem}[Lower unit window]\label{lw:thm:lw}
There are $a_{\rm LW}>0$ and $T_0$ such that, for every $T\ge T_0$
and every $u/\tilde t\in J$,
\begin{equation}\label{lw:eq:lw}
 \boxed{\displaystyle
 \PT\bigl(\cC_L\cap\{S_L\in[u,u+1]\}\bigr)
 \ge a_{\rm LW}\,t^{-1/2}e^{-u^2/\tilde t}.}
 \tag{LW}
\end{equation}
The constants are uniform in $u$, in $1\le q\le T^{\theta_0}$ and
$\chi\bmod q$, in $\eps\in\{-1,1\}$, and in $x_L\in[x_*,x_*+1)$. The
threshold depends only on $\theta_0,J,R,K_0,\Ac,\beta_*$ and the
common Taylor calibration; $M$ enters through that calibration. One
may take
\[
 a_{\rm LW}=\frac1{80}e^{-R-C_R},
\]
where $C_R$ is any constant, independent of $q$ and $\chi$, with
$\log\E e^{\lambda S_L^*}\ge\lambda^2\tilde t/4-C_R$ for
$|\lambda|\le R$. No hypothesis on zeros is used.
\end{theorem}

The parameters are fixed for each compact $J$ before $T$ tends to
infinity. The terminal shift ranges throughout its full interval.
The proof occupies the next four sections.

\section{Eightfold Taylor inflation}\label{lw:sec:inflation}

Replace each $p^{-i\tau}$ by an independent Haar variable $U_p$.
Expectation $\E$ refers to this model; stars are suppressed during model
calculations. Set
\[
 \lambda=2u/\tilde t,\quad v=u+\tfrac12,\quad
 Z_i=\E e^{\lambda Y_i},\quad Z=\prod_iZ_i,
 \quad d\Pp_\lambda=Z^{-1}e^{\lambda S_L}\,d\Pp.
\]
This tilt is still independent over primes and blocks. Write
$E_m(z)=\sum_{j=0}^m z^j/j!$ and set
\begin{equation}\label{lw:eq:polynomials}
 m_i=8N_i,\quad z=\lambda/2+i\xi,\quad
 P_i(\xi)=E_{m_i}(zY_i),\quad Q_i(\xi)=e^{zY_i},\quad
 D_i(\xi)=Y_i\frac{(zY_i)^{m_i}}{m_i!}.
\end{equation}
Only the polynomial degrees are multiplied by eight. The constants
$\Ac$, $K_0$ and $\beta_*$ remain those fixed in
Section~\ref{sec:parameters}.

\begin{lemma}[Inflation]\label{lw:lem:inflation}
Put $\epsilon_i=e^{-N_i}$. Uniformly for $|\lambda|\le R$ and
$|\xi|\le4$,
\begin{align}
 |P_i(\xi)e^{-zY_i}-1|&\le\epsilon_i
                    &&\text{if }|Y_i|\le B_i,\label{lw:eq:pointwise-taylor}\\
 \|P_i-Q_i\|_2,\quad \|Y_i(P_i-Q_i)\|_2,\quad\|D_i\|_2
                    &\le\epsilon_i\sqrt{Z_i}.\label{lw:eq:l2-taylor}
\end{align}
The norms are under the untilted model. The choice of $\Ac$ in
Lemma~\ref{lem:taylor} also gives
\begin{equation}\label{lw:eq:ac-lower}
 \Ac\ge eK_0\rho_0,\qquad \rho_0=\sqrt{R^2/4+1}.
\end{equation}
\end{lemma}

\begin{proof}
We use the positive majorants fixed in Lemma~\ref{lem:taylor},
including \eqref{eq:calibrated-moment} and
\eqref{eq:calibrated-cap}. The majorants themselves, rather than
only the actual Taylor error, provide the needed degree comparison.
For $j\ge1$ put
\[
 \mathcal R_i(j,\rho)=\frac{\rho^j}{j!}
 \left\{Z_i^{-1/2}(\E|Y_i|^{2j})^{1/2}
                  +(\E_\lambda|Y_i|^{2j})^{1/2}\right\}.
\]
The bounded-variable moment estimate used in Lemma~\ref{lem:taylor},
with the same calibrated constant $C$, gives
\begin{equation}\label{lw:eq:majorant}
 \mathcal R_i(j,\rho)\le
 e^{E_i}\left[C\frac{\rho}{\rho_0}
          \frac{\nu_i+\sqrt{j\nu_i}}j\right]^j
 =:\mathcal B_i(j,\rho),\qquad E_i=C\nu_i\ge0.
\end{equation}
With $k_0=N_i+1$, the already chosen $\Ac$ ensures
\begin{equation}\label{lw:eq:calibration}
 \mathcal B_i(k_0,\rho_0)\le e^{-N_i},\qquad
 (1+e^{RB_i/2})\frac{(\rho_0B_i)^{k_0}}{k_0!}\le e^{-N_i}.
\end{equation}
The first bound is precisely the unweighted specialization
\eqref{eq:calibrated-moment}; the second is
\eqref{eq:calibrated-cap}.

Let $N=N_i$, $\rho_4=\sqrt{R^2/4+16}\le4\rho_0$ and
$j\ge8N+1$. The bracket in \eqref{lw:eq:majorant} is at most
$4\sqrt{k_0/j}$ times its value at $(k_0,\rho_0)$. Hence
\[
 \log\mathcal B_i(j,\rho_4)
 \le j\left\{\log(4\sqrt{k_0/j})-\frac N{k_0}\right\}
       -(j/k_0-1)E_i\le-j/4.
\]
Indeed $N\ge401$, so in particular $j\ge6k_0$, $N/k_0\ge3/4$,
and $\log(4/\sqrt6)<1/2$. The integral Taylor remainder now gives,
with $m=8N$, $k=m+1$,
\begin{align*}
 \|P_i-Q_i\|_2/\sqrt{Z_i}&\le e^{-k/4},\\
 \|Y_i(P_i-Q_i)\|_2/\sqrt{Z_i}
      &\le(k+1)\rho_4^{-1}e^{-(k+1)/4},\\
 \|D_i\|_2/\sqrt{Z_i}&\le k\rho_4^{-1}e^{-k/4}.
\end{align*}
The factorial multipliers in the last two lines are exact. As
$\rho_4\ge4$ and $xe^{-x}\le e^{-x/2}$ for $x\ge2$, all three
right sides are at most $e^{-N}$.

For the pointwise assertion, let $A_i=1+e^{RB_i/2}$ and
$x=\rho_0B_i$. Equation \eqref{lw:eq:calibration} and
$k_0!\le k_0^{k_0}$ give
$x\le k_0e^{-(N+\log A_i)/k_0}$. By $k!\ge(k/e)^k$,
\[
 \log\left(A_i\frac{(\rho_4B_i)^k}{k!}\right)
 \le k\left\{\log(4ek_0/k)-N/k_0\right\}
       -(k/k_0-1)\log A_i
 \le-k/8\le-N.
\]
Here $\log(2e/3)-3/4<-1/8$. This proves the relative cap bound
for both signs of the real part of $zY_i$.

In the explicit factorial calibration of Lemma~\ref{lem:taylor},
the exponential multiplying $2e\rho_0K_0/\Ac$ is at least one.
Hence $\Ac\ge2e^2\rho_0K_0$, which implies \eqref{lw:eq:ac-lower}.
\end{proof}

We will use two quantitative consequences. First,
\begin{equation}\label{lw:eq:epsilon}
 \delta=\sum_{i\le L}\epsilon_i
 \le\frac{e^{-401}}{1-e^{-200}}<2e^{-400},\qquad
 \rho=\delta+\sqrt t\,\epsilon_1\le3e^{-400}
\end{equation}
after increasing $T_0$ only. The geometric bound follows from
$c_L\ge400$ and $e^{j/2}\ge1+j/2$. Second, by
$1+K^2\le\rho_0^2$, \eqref{lw:eq:cutoff} and \eqref{lw:eq:ac-lower},
\begin{equation}\label{lw:eq:cap-smallness}
 E_*:=\frac{(1+K^2)\beta_*}{K_0^2}
 \le\frac1{e^2(32\Cg)^2K_0^4},\qquad
 \frac{4E_*}{1-e^{-1}}<\frac1{40},\quad
 \frac{240E_*}{1-e^{-1}}<1.
\end{equation}
For the last inequalities even $\Cg\ge1$, $e^2\ge1$, $K_0\ge4$
and $1-e^{-1}>1/2$ suffice. These comparisons are uniform over the
permitted slope range, with the calibrated $\Ac$ fixed throughout.

\section{The tilted local limit and the cost of the caps}\label{lw:sec:local}

Mertens' formula \eqref{lw:eq:mertens} with a constant, through
Lemma~\ref{lem:clock}(iii), yields the precise later-block second
moments below.

Write $y_p=\eps\Ima(\chi(p)U_p)/\sqrt p$ for $p\nmid q$, omitting
the primes dividing $q$. At a retained prime, rotation removes
$\eps\chi(p)$, and the partition function is $I_0(\lambda/\sqrt p)$, with
$I_0(s)=(2\pi)^{-1}\int_0^{2\pi}e^{s\sin\theta}\,d\theta$.
Uniformly on $|\lambda|\le R$,
\begin{equation}\label{lw:eq:prime-cumulants}
 \log I_0(\lambda/\sqrt p)=\lambda^2/(4p)+O_R(p^{-2}),\quad
 \E_\lambda y_p=\lambda/(2p)+O_R(p^{-2}),\quad
 \Var_\lambda y_p=1/(2p)+O_R(p^{-2}).
\end{equation}
Consequently
\begin{equation}\label{lw:eq:cumulants}
 \log Z=\lambda^2\tilde t/4+O_R(1),\quad
 \mu:=\E_\lambda S_L=\lambda\tilde t/2+O_R(1),\quad
 \sigma^2:=\Var_\lambda S_L=\tilde t/2+O_R(1).
\end{equation}
The sums run over $p\le T_L$ with $p\nmid q$, and the implied
constants do not depend on $q$ or $\chi$. For example, if
\[
 b_R=\sup_{p,\ |\lambda|\le R}
       p^2\left|\log I_0(\lambda/\sqrt p)-\lambda^2/(4p)\right|<\infty,
\]
then an admissible constant in Theorem~\ref{lw:thm:lw} is
\begin{equation}\label{lw:eq:mgf-constant}
 C_R=b_R\sum_p p^{-2}.
\end{equation}

We need an absolute damping constant. Delete the fixed primes at most
$\max(16R^2,4096)$ and those dividing $q$. At a remaining prime the
tilted density relative to Haar measure is at least $e^{-1/2}>1/2$,
so $\Var_\lambda y_p\ge1/(4p)$. For $|\tau|\le8$, independent
tilted copies $y_p,y_p'$ satisfy $|\tau(y_p-y_p')|\le1$. Therefore
\[
 1-|\E_\lambda e^{i\tau y_p}|^2
 \ge\tfrac13\tau^2\E_\lambda(y_p-y_p')^2
 \ge\tau^2/(6p).
\]
The retained first-block reciprocal mass is at least $t_1/2$ for
large $T$. Multiplication gives the convenient weaker bound
\begin{equation}\label{lw:eq:damping}
 |\E_\lambda e^{i\tau Y_1}|\le e^{-t_1\tau^2/32}
                  \qquad(|\tau|\le8).
\end{equation}
The constant $1/32$ is absolute; by Lemma~\ref{lem:clock}(i) the
threshold depends only on $R$ and $\theta_0$.

\begin{lemma}[Uniform Fourier local limit]\label{lw:lem:local}
Use the convention $\widehat f(\xi)=(2\pi)^{-1}\int f(x)e^{-i\xi x}\,dx$.
If $f\in L^1(\mathbb R)$ has continuous Fourier transform supported in
$[-8,8]$, then
\begin{equation}\label{lw:eq:local}
 \sqrt t\,\E_\lambda f(S_L-v)
       \longrightarrow\frac1{\sqrt\pi}\int_{\mathbb R}f(x)\,dx.
\end{equation}
For every later block, uniformly also in $2\le i\le L$,
\begin{equation}\label{lw:eq:insert}
 \sqrt t\,\E_\lambda[Y_i^2 f(S_L-v)]
 =\frac{\kappa^2+1/2}{\sqrt\pi}\int f+o(1),\qquad\kappa=u/\tilde t.
\end{equation}
All limits are uniform for $\kappa\in J$, for $1\le q\le T^{\theta_0}$
and $\chi\bmod q$, for both signs and over the full terminal-shift
interval.
\end{lemma}

\begin{proof}
The centered one-prime third absolute moments have bounded sum,
since each is at most $8p^{-3/2}$. For sufficiently small fixed $|\tau|$,
expanding the centered characteristic functions and their logarithms
gives
\[
 \log\E_\lambda e^{i\tau(S_L-\mu)}
       =-\sigma^2\tau^2/2+O_R(|\tau|^3).
\]
Since $\mu-v=O_R(1)$ and $\sigma^2/t\to1/2$ uniformly
(Lemma~\ref{lem:clock}(iii)), putting $\tau=z/\sqrt t$ gives
\[
 e^{-izv/\sqrt t}\E_\lambda e^{izS_L/\sqrt t}
                        \longrightarrow e^{-z^2/4}
\]
uniformly on bounded $z$ intervals. Fourier inversion gives the continuous
representative of $f$. The first block contains a retained prime for
large $T$ by Lemma~\ref{lem:clock}(i); its tilted sine has a density,
so changing the representative leaves either expectation unchanged.
Fourier inversion therefore expresses the
left side of \eqref{lw:eq:local} as the integral of this expression times
$\widehat f(z/\sqrt t)$ over $|z|\le8\sqrt t$. By
\eqref{lw:eq:damping} and $t_1\ge t/2$, its absolute value is bounded by
$\|\widehat f\|_\infty e^{-z^2/64}$. Dominated convergence, uniformly
after splitting into bounded and unbounded $z$, gives
$2\sqrt\pi\widehat f(0)=\pi^{-1/2}\int f$.

For $i\ge2$, Lemma~\ref{lem:clock}(iii) gives
$\sum_{p\in\mathcal P_i,\,p\nmid q}p^{-1}=\tilde t_i-\tilde t_{i-1}=1+o(1)$
uniformly. Equations
\eqref{lw:eq:prime-cumulants} and the bounded-variable exponential
estimate used in Lemma~\ref{lem:taylor} give
\[
 \E_\lambda Y_i^2=\kappa^2+1/2+o(1),\qquad
 \E_\lambda|Y_i|^3=O_R(1).
\]
Inserting $Y_i^2$ in the Fourier integral changes its block factor at
$\tau=z/\sqrt t$ by a ratio
$\E_\lambda Y_i^2+O_R(|\tau|)$ on bounded $z$ intervals. The
unmodified block factor is then bounded away from zero. For domination
on the full frequency interval, the inserted factor is bounded by
$\E_\lambda Y_i^2=O_R(1)$, while the untouched first block still
supplies \eqref{lw:eq:damping}. The same integral argument proves
\eqref{lw:eq:insert}. All estimates are uniform in the stated parameters.
\end{proof}

Set
\begin{equation}\label{lw:eq:kernel}
 h(\xi)=\frac{(4-|\xi|)_+}{16},\quad
 k(x)=\int h(\xi)e^{i\xi x}\,d\xi
      =\left(\frac{\sin2x}{2x}\right)^2,\quad
 g(x)=(1-4x^2)k(x)^2.
\end{equation}
Here $k(0)=1$. The transforms of $k^2$ and $x^2k^2$ are respectively
$h*h$ and $-(h*h)''=-h'*h'$, continuous and supported in $[-8,8]$.
Both functions are integrable. Plancherel yields
\begin{equation}\label{lw:eq:kernel-masses}
 \int h^2=\frac16,\quad \int|h'|^2=\frac1{32},\quad
 \int k^2=\frac\pi3,\quad \int x^2k^2=\frac\pi{16},\quad
 \int g=\frac\pi{12}>0.
\end{equation}
Write
\[
 F_v^\circ=e^{\lambda S_L/2}k(S_L-v),\quad
 U=Z/\sqrt t,\quad G_0=\sqrt\pi/12.
\]
The local limit and insertion now say
\begin{align}
 \E\{|F_v^\circ|^2-4|(S_L-v)F_v^\circ|^2\}
                           &=(G_0+o(1))U,\label{lw:eq:core-exact}\\
 \E[Y_i^2|F_v^\circ|^2]
             &=\left\{\frac{\sqrt\pi}{3}(\kappa^2+\tfrac12)+o(1)\right\}U
                                      \quad(i\ge2).\label{lw:eq:cap-exact}
\end{align}
Since $B_i^{-2}=\beta_i/K_0^2$ and
$\sum_{i=2}^L\beta_i\le\beta_* /(1-e^{-1})$, the ratio of the
leading later-cap penalty to $G_0U$ is at most
\begin{equation}\label{lw:eq:cap-ratio}
 \frac{4(\kappa^2+1/2)}{K_0^2}\sum_{i=2}^L\beta_i
 \le\frac{4E_*}{1-e^{-1}}<\frac1{40}.
\end{equation}
The uniform $o(1)$ errors remain $o(U)$ on summing with these weights.

\section{A signed lower certificate and its localized errors}\label{lw:sec:certificate}

Define the finite Laurent polynomials
\begin{equation}\label{lw:eq:certificate}
 F_v=\int_{-4}^4 h(\xi)e^{-i\xi v}\prod_{i\le L}P_i(\xi)\,d\xi,
 \qquad
 H_v=\left(1-4(S_L-v)^2-\sum_{i\le L}(Y_i/B_i)^2\right)|F_v|^2.
\end{equation}
The integral changes only coefficients. If the window or any cap
fails, the real bracket is nonpositive. On the capped window,
\eqref{lw:eq:pointwise-taylor}, $\int h=1$ and \eqref{lw:eq:epsilon} give
$|F_v|\le2e^{\lambda S_L/2}$. Thus, globally on the arithmetic interval,
\begin{equation}\label{lw:eq:global-minorant}
 \boxed{\displaystyle
 \one_{\cC_L}\one_{\{S_L\in[u,u+1]\}}
       \ge\frac{e^{-\lambda u}}{4e^R}H_v.}
\end{equation}
No approximation is needed outside the caps.

We next preserve the model positivity at scale $U$, rather than at the
larger scale $Z$. The details below give absolute constants, so a fixed
Taylor error is never compared to an unspecified large constant.

For $f=h$ or $h'$, let $F_f$ be the integral defining $F_v$ with $h$
replaced by $f$, and let $F_f^\circ$ use $Q_i$ instead of $P_i$. Define
normalized block kernels
\[
 K_i(\xi,\zeta)=Z_i^{-1}\E[P_i(\xi)\overline{P_i(\zeta)}],\qquad
 K_i^\circ(\xi,\zeta)=\E_\lambda e^{i(\xi-\zeta)Y_i}.
\]
By \eqref{lw:eq:l2-taylor},
$|K_i-K_i^\circ|\le3\epsilon_i$ and $|K_i|\le1+3\epsilon_i$.
Products of later kernels have magnitude at most $2$ and differ from
their exact products by at most $4\delta$. Consequently, with
$d(\tau)=e^{-t_1\tau^2/32}$,
\begin{equation}\label{lw:eq:kernel-error}
 \left|\prod_i K_i-\prod_iK_i^\circ\right|
                         \le6\epsilon_1+4\delta\,d(\xi-\zeta).
\end{equation}
The elementary convolution bounds
\begin{equation}\label{lw:eq:j-bounds}
 \iint|h(\xi)h(\zeta)|d(\xi-\zeta)\,d\xi d\zeta\le4/\sqrt t,
 \quad
 \iint|h'(\xi)h'(\zeta)|d(\xi-\zeta)\,d\xi d\zeta\le1/\sqrt t
\end{equation}
follow from $\|h\|_1=1$, $\|h\|_\infty=1/4$,
$\|h'\|_1=1/2$, $\|h'\|_\infty=1/16$, and $t_1\ge t/2$.
Equations \eqref{lw:eq:kernel-error}--\eqref{lw:eq:j-bounds} give
\begin{equation}\label{lw:eq:f-errors}
 \big|\E|F_h|^2-\E|F_h^\circ|^2\big|\le22\rho U,
 \qquad
 \big|\E|F_{h'}|^2-\E|F_{h'}^\circ|^2\big|\le6\rho U.
\end{equation}

Multiplying an unlocalized Taylor error by $S_L-v$ would not suffice.
Instead use $P_i'=iY_iP_i-iD_i$. Integration by parts, with $h$ zero
at its endpoints and absolutely continuous, gives the exact identity
\begin{equation}\label{lw:eq:ibp}
 (S_L-v)F_h=iF_{h'}+\sum_i R_i,\qquad
 R_i=\int h(\xi)e^{-i\xi v}D_i(\xi)\prod_{j\ne i}P_j(\xi)\,d\xi.
\end{equation}
For the exact exponentials, $(S_L-v)F_h^\circ=iF_{h'}^\circ$.
For $i\ge2$, the first block remains unmodified in the norm of $R_i$.
Using \eqref{lw:eq:l2-taylor} and \eqref{lw:eq:j-bounds},
\[
 \|R_i\|_2^2\le(8+6\sqrt t\epsilon_1)\epsilon_i^2U
                    \le9\epsilon_i^2U.
\]
For the first block, direct Minkowski gives
$\|R_1\|_2\le2\epsilon_1\sqrt Z$. Therefore
\begin{equation}\label{lw:eq:defect}
 \Big\|\sum_iR_i\Big\|_2\le3\rho\sqrt U.
\end{equation}
By \eqref{lw:eq:damping}--\eqref{lw:eq:j-bounds},
$\|F_{h'}^\circ\|_2^2\le U$; by \eqref{lw:eq:f-errors},
$\|F_{h'}\|_2\le2\sqrt U$. Since $\rho<1/9$,
\[
 \left|\E|(S_L-v)F_h|^2-\E|(S_L-v)F_h^\circ|^2\right|
 \le(6\rho+12\rho+9\rho^2)U\le19\rho U.
\]
Together with \eqref{lw:eq:f-errors}, this proves
\begin{equation}\label{lw:eq:core-error}
 \left|\E\{|F_h|^2-4|(S_L-v)F_h|^2\}
       -\E\{|F_h^\circ|^2-4|(S_L-v)F_h^\circ|^2\}\right|
                                      \le100\rho U.
\end{equation}

For completeness the cap insertion has an equally uniform error.
Let $M_c=K+2$, distinct from the overlap parameter $M$.
For $i\ge2$ and large $T$, the inserted block has
\[
 \|Y_iQ_i\|_2\le(K+1)\sqrt{Z_i},\quad
 \|Y_iP_i\|_2\le M_c\sqrt{Z_i},\quad
 \|Y_i(P_i-Q_i)\|_2\le\epsilon_i\sqrt{Z_i}.
\]
Thus its normalized kernel has magnitude at most $M_c^2$, and its
error is at most $2M_c\epsilon_i$. A telescoping product gives a
later-product error at most $6M_c^2\delta$ and magnitude at most
$2M_c^2$. The same first-block separation and the $h$ bound in
\eqref{lw:eq:j-bounds} give
\begin{equation}\label{lw:eq:cap-error}
 \left|\E[Y_i^2|F_h|^2]-\E[Y_i^2|F_h^\circ|^2]\right|
                      \le30M_c^2\rho U\qquad(i\ge2).
\end{equation}
After division by $B_i^2$ and summation this is at most $\rho U$:
indeed $(K+2)^2\le8(1+K^2)$ and \eqref{lw:eq:cap-smallness} apply.

The first cap is treated directly. The same moment estimate gives
$\|Y_1P_1\|_2\le C_R(t_1+1)\sqrt{Z_1}$. Minkowski and independence
then give
\begin{equation}\label{lw:eq:first-cap}
 B_1^{-2}\E[Y_1^2|F_h|^2]
 \le C_RK_0^{-2}(t_1+1)^2e^{-t_1}Z=o(U).
\end{equation}
Combining \eqref{lw:eq:core-exact}, \eqref{lw:eq:cap-ratio} and
\eqref{lw:eq:core-error}--\eqref{lw:eq:first-cap} now proves
\begin{equation}\label{lw:eq:positive}
 \E H_v\ge\{(39/40)G_0-101\rho+o(1)\}U\ge U/10
\end{equation}
for all sufficiently large $T$. Here $(39/40)G_0>0.14$ and
$\rho\le3e^{-400}$, with all numerical error constants displayed.
The sum of the absolute square terms also satisfies
\begin{equation}\label{lw:eq:absolute-mass}
 \E|F_v|^2+4\E|(S_L-v)F_v|^2
               +\sum_iB_i^{-2}\E|Y_iF_v|^2\le10U
\end{equation}
eventually, by the same estimates. No factor proportional to $L$ occurs.

\section{Arithmetic transfer for the lower unit window}\label{lw:sec:transfer}

Lemma~\ref{lem:mv}, based on Montgomery--Vaughan's separated-frequency
Hilbert inequality \cite{MV1974}, says that a Laurent
polynomial of symmetric rational height at most $T^h$ satisfies
\begin{equation}\label{lw:eq:mv}
 |\ET|F|^2-\E|F^*|^2|\le C_{\rm MV}T^{2h-1}\E|F^*|^2.
\end{equation}
It is uniform in all coefficients and degrees, after combining equal
rational frequencies. This is valid for the character phases used here.

Put $\mathcal N=\sum_{i\le L}\beta_iN_i$. The height of $F_v$ is
at most $T^{8\mathcal N}$, and multiplication by $S_L-v$ or any
$Y_i$ costs at most an extra $T^{\beta_L}$. The full exponent in
the squared transfer is consequently at most
\begin{align}
 16\mathcal N+2\beta_L
 &\le16\Ac\beta_1t_1+16\Cg\Ac\sqrt{\beta_L}+2\beta_L\notag\\
 &\le16\Ac\beta_1t_1+\tfrac12+2\beta_*\le\tfrac34.
 \label{lw:eq:height}
\end{align}
The last step only enlarges $T_0$ so that
$16\Ac\beta_1t_1\le1/8$, using $\beta_*\le400^{-2}$.
The first inequality is the explicit ceiling estimate
\eqref{eq:certificate-cost}. Both numerator and denominator
frequencies are counted. The factor-eight construction uses the
explicit $(32\Cg\Ac)^{-2}$ cutoff in
\eqref{lw:eq:cutoff}.

Apply \eqref{lw:eq:mv} separately to the nonnegative squares comprising
$H_v$, and add their absolute errors. Equations
\eqref{lw:eq:absolute-mass}--\eqref{lw:eq:height} imply
\[
 |\ET H_v-\E H_v^*|\le10C_{\rm MV}T^{-1/4}U.
\]
Enlarge $T_0$ once more. By \eqref{lw:eq:positive},
$\ET H_v\ge U/20$. The global minorant \eqref{lw:eq:global-minorant}
therefore gives
\[
 \PT(\cC_L,S_L\in[u,u+1])
 \ge\frac{e^{-R}}{80}e^{-\lambda u}\frac Z{\sqrt t}
 \ge\frac{e^{-R-C_R}}{80}\frac{e^{-u^2/\tilde t}}{\sqrt t},
\]
since $\lambda=2u/\tilde t$. This proves Theorem~\ref{lw:thm:lw}.
\hfill$\square$

\section{Fixed-width windows for the argument}\label{sec:local-window}

\begin{proof}[Proof of Theorem~\ref{thm:local-window}]
Write $Q=\cL_T$, $a=|v|$ and $\kappa=a/Q$. Choose $\eps$ so that
$\eps v=a$, taking $\eps=1$ when $v=0$. Then
$|X_\chi-v|=|\eps X_\chi-a|$. Use the signed walk \eqref{eq:reference}
for this $\eps$ throughout. Apply Sections~\ref{sec:parameters}
and~\ref{lw:sec:statement} with slope endpoint $K+1$, so that
$J=[0,K+1]$ and $R\ge2(K+1)$ before the degrees and cutoff are fixed.
All constants depend only on $K,\theta_0$ and the data of \textup{(OS)}.
Put
\[
 d(a)=Q^{-1/2}\exp(-a^2/Q).
\]

\emph{The upper bound.}
Fix $w>0$ and let $E=\{|\eps X_\chi-a|\le w\}$. At any level $j$,
the unit-window factor in Proposition~\ref{prop:comparison} at
$u$ and tilt $\lambda=2\kappa$ satisfies
\begin{equation}\label{eq:local-fixed-tilt}
 \frac{e^{\kappa^2\tilde t_j-2\kappa u}}{\sqrt{t_j}}
 =d(a)\sqrt{\frac{Q}{t_j}}\,
   e^{-\kappa^2(Q-\tilde t_j)+2\kappa(a-u)}
 \ll d(a)e^{2K|a-u|}.
\end{equation}
Indeed $\tilde t_j\le Q$ and $t_j\ge n/2$, whereas $Q\le n+1$.
This includes $\kappa=0$.

Start with $G_1=\{|S_1|\le c_1/2\}$. For $2\le\ell\le L$ put
\[
 A_\ell=\{|Y_\ell|\le2c_\ell\},\qquad
 G_\ell=G_{\ell-1}\cap A_\ell\cap\{|S_\ell-a|\le w+c_\ell\}.
\]
Since $a\le K(n+1)\le c_1/2$ for large $T$, on $G_j$ we have
$|S_j-a|\le w+c_j$ at every level, including $j=1$. The increment
caps \eqref{eq:caps} hold through $j$ on $G_j$, and through $\ell$
on $G_{\ell-1}\cap A_\ell$. The initial error \eqref{eq:initial-error}
is $O(d(a))$, since $d(a)\ge(n+1)^{-1/2}e^{-K^2(n+1)}$.

Consider first $G_{\ell-1}\cap A_\ell^c$. Cover the possible values of
$S_{\ell-1}$ by the integer unit intervals meeting
$[a-w-c_{\ell-1},a+w+c_{\ell-1}]$. There are $O(1+w+c_{\ell-1})$
such intervals, each with $|a-u|\le w+c_{\ell-1}+1$ at its left
endpoint. Use the weighted comparison at level $\ell-1$ with
$r=r_\ell=\lfloor c_\ell/D\rfloor$ and $W=Y_\ell$. Its proxy is
$O(1)$ and its absolute overlap is zero. Markov's inequality divides
the $2r$th moment by $(2c_\ell)^{2r}$. Since
$c_{\ell-1}=e^{1/2}c_\ell$, \eqref{eq:local-fixed-tilt} gives
\begin{equation}\label{eq:local-increment-exit}
 \PT(G_{\ell-1}\cap A_\ell^c)
 \ll_{K,w}d(a)(1+c_\ell)e^{C_Kc_\ell}
              \left(\frac{Cr_\ell}{c_\ell^2}\right)^{r_\ell}.
\end{equation}

For the other exits use $W_\ell^\pm=P^\pm_{y_\ell}-S_\ell$ with
$y_\ell=T^{A/c_\ell}$. The two signs of \textup{(OS)} in
\eqref{eq:majorant-final} give
\[
 S_\ell+W_\ell^- -Bc_\ell/A\le \eps X_\chi
                 \le S_\ell+W_\ell^+ +Bc_\ell/A.
\]
On $E$, if $S_\ell<a-w-c_\ell$ then $W_\ell^+>7c_\ell/8$.
If $S_\ell>a+w+c_\ell$ then $-W_\ell^->7c_\ell/8$. Here $A\ge8B$.
In both cases the normalization of the two weights at zero is
$-i\eps$, so Lemma~\ref{lem:coefficients} applies to each correction
from the same $S_\ell$. On $G_{\ell-1}\cap A_\ell$ the triangle
inequality gives
\[
 |S_\ell-a|\le w+c_{\ell-1}+2c_\ell\le w+4c_\ell .
\]
Cover this finite interval as above. On each capped unit window use
Proposition~\ref{prop:comparison} separately for $W_\ell^+$ and
$W_\ell^-$, each with proxy $\Sigma_\ell=C_3M(x_\ell)$ from
\eqref{eq:merged-proxy} and moment order $2r_\ell$. Dividing by
$(7c_\ell/8)^{2r_\ell}$ and summing the two bounds gives
\begin{align}
 &\PT\bigl(E\cap G_{\ell-1}\cap A_\ell
                    \cap\{|S_\ell-a|>w+c_\ell\}\bigr)\notag\\
 &\qquad\ll_{K,w} d(a)(1+c_\ell)e^{C_Kc_\ell}
          \left(\frac{Cr_\ell M(x_\ell)}{c_\ell^2}\right)^{r_\ell}.
 \label{eq:local-barrier-exit}
\end{align}
For each comparison the certificate cost $2\sum\beta_iN_i$ through
its level is below $1/8$ by \eqref{eq:certificate-cost}.
The additional height cost is at most
$2r_\ell\beta_\ell\le2\sqrt{\beta_\ell}/D<1/8$ for $Y_\ell$
and $2r_\ell A/c_\ell\le 2A/D<1/8$ for either correction.
Also $1\le r_\ell\le c_\ell\le c_{\ell-1}$ and the overlaps satisfy
\eqref{eq:absolute-overlap}. Thus \eqref{eq:height-budget} holds
for each use separately.

After division by $d(a)$ the two exit bounds are summable uniformly in
$T$. Indeed $c_\ell=e^{x_\ell/2}$,
$c_\ell/(2D)\le r_\ell\le c_\ell/D$, and
$M(x)\ll_A 1+x$. For sufficiently large $x=x_\ell$ the logarithm
of the common upper bound is at most
\[
 O_{K,w}(1+c_\ell)
 +r_\ell\log\!\left(\frac{C_A(1+x)}{c_\ell}\right)
 \le -\frac{x e^{x/2}}{16D}.
\]
The $x_\ell=x_L+L-\ell$ have spacing one, with
$x_L\in[x_*,x_*+1)$ by \eqref{lw:eq:ladder}. The bounded range of $x_\ell$
contributes a fixed constant as well. Finally, on $G_L$ the interval
$[a-w-c_L,a+w+c_L]$ is covered by $O_{w,\beta_*}(1)$ unit windows
whose left endpoints have $|a-u|\le w+c_L+1$. The unweighted
comparison and \eqref{eq:local-fixed-tilt} give $\PT(G_L)\ll d(a)$.
Every point of $E$ is in $G_1^c$, has a first exit at some
$2\le\ell\le L$, or belongs to $G_L$. This proves the upper bound.

\emph{The lower bound.}
At the terminal level put $\tau=\tilde t_L$ and $t=t_L$.
By Lemma~\ref{lem:clock}(iii) and \eqref{lw:eq:ladder}, $0\le Q-\tau\le x_*+2$ and
$t\asymp\tau\asymp Q$. As $Q\to\infty$ uniformly in $q$ by
Lemma~\ref{lem:clock}, we have $a/\tau\in[0,K+1]$ for all large $T$.
On the same signed walk take the single capped window
$U=\cC_L\cap\{S_L\in[a,a+1]\}$.
The signed certificate in Theorem~\ref{lw:thm:lw} gives
\[
 \PT(U)\ge a_{\rm LW}\,t^{-1/2}e^{-a^2/\tau}.
\]
Both terminal corrections have proxy
$\Sigma_*=C_3(1+\log A+(x_*+1)/2)$. With
$\lambda=2a/\tau$, $r=1$, the comparison on $U$ and Markov's
inequality give
\[
 \PT\bigl(U\cap\{\max(|W_L^-|,|W_L^+|)>H\}\bigr)
 \le 2C_{\rm cp}\Sigma_*H^{-2}\,t^{-1/2}e^{-a^2/\tau}.
\]
Here the factor two is a union bound. The absolute overlap and height
estimates already verified apply to each weight with $r=1$.
Choose $H\ge1$ after $a_{\rm LW}$ so that
$2C_{\rm cp}\Sigma_* H^{-2}\le a_{\rm LW}/2$. On the remaining
event the envelopes put $\eps X_\chi$ in
$[a-H-d_*,a+1+H+d_*]$, where $d_*=(B/A)e^{(x_*+1)/2}$.
Taking $w_0=2+H+d_*$ therefore gives, for all $w\ge w_0$,
\[
 \PT\{|\eps X_\chi-a|<w\}\ge
 \tfrac12 a_{\rm LW}\,t^{-1/2} e^{-a^2/\tau}\gg d(a).
\]
The last comparison follows from
$0\le a^2/\tau-a^2/Q\le K^2(Q/\tau)(x_*+2)$.
All choices are uniform in $v,q,\chi$ and the sign, including the
zero tilt at $v=0$.
\end{proof}

\section{From windows to tails and moments}\label{lw:sec:tails}

We first bound the contribution of heights beyond any sufficiently large
fixed slope. This estimate uses only \textup{(OS)} and mean-square
transfer, so it applies before the sharp tail bounds are established.

\begin{lemma}[Far tail]\label{lem:far-tail}
Write $Y=\eps X_\chi$ and $\cL=\cL_T$. Under the hypotheses of
Theorem~\ref{thm:uniform}, for each $k>0$ and every fixed
\[
 K\ge K_0'(k):=\max\Bigl\{k+2,\;6C_w,\;
      \frac{3C_w^2}{C_2}e^{8C_2(1+\theta_0)(k+1)}\Bigr\},
\]
we have $\ET[e^{2kY}\one_{\{Y>K\cL\}}]\ll e^{-2K\cL}$ for large $T$,
uniformly in $q,\chi,\eps$.
\end{lemma}
\begin{proof}
Take the fixed cutoff $y_1=\max\{y_0,2\}$. Bounding its prime and
prime-square terms individually gives, outside a null set,
\[
 Y\le V_{\max}:=C_2\frac{\log(qT)}{\log y_1}+C_1,
 \qquad C_1=C_0+C_w\sum_{p\le y_1}(p^{-1/2}+p^{-1}).
\]
Let $V=K\cL+j$ with an integer $j\ge0$ and $V\le V_{\max}$; there are
finitely many such $V$. Define $y_V$ by
$\log y_V=2C_2\log(qT)/V$. For large $T$ we have $y_1\le y_V\le T$,
because $V_{\max}<2C_2\log(qT)/\log y_1$ and
$V\ge K\cL>2C_2(1+\theta_0)\ge2C_2\log(qT)/\log T$. Hence
\eqref{eq:majorant} at the truncation $y_V$ gives, outside a null set,
\[
 Y\le P_{\eps,y_V}+\frac V2+C_0 .
\]
Write $P_{\eps,y_V}=\Rea A+\Rea B$, where
$A(t)=\sum_{p\le y_V}a_pp^{-it}$ with
$a_p=\chi(p)p^{-1/2}w_\eps(\log p/\log y_V)$, and $B$ is the sum over
prime squares. Since $y_V\le T$ and $\chi(p)=0$ for $p\mid q$,
\[
 |B|\le\frac{C_w}2\sum_{\substack{p\le T\\p\nmid q}}\frac1p
 =\frac{C_w}2\cL,
 \qquad \sum_{p\le y_V}|a_p|^2\le C_w^2\cL.
\]
Since $V\ge K\cL\ge6C_w\cL$, for large $T$ we have
$V/2-C_0-|B|\ge V/3$. Therefore, up to a null set,
$\{Y>V\}\subset\{|A|>V/3\}$.

Let $m=\lfloor V/(4C_2(1+\theta_0))\rfloor$, which is at least $1$ for
large $T$. Then $A^m=\sum_{N\le y_V^m}c_NN^{-it}$ and
$y_V^m\le(qT)^{1/(2(1+\theta_0))}\le T^{1/2}$. If
$N=\prod_ip_i^{\alpha_i}$ with $\sum_i\alpha_i=m$, unique factorization
gives $c_N=(m!/\prod_i\alpha_i!)\prod_ia_{p_i}^{\alpha_i}$, so
$|c_N|^2\le m!\,(m!/\prod_i\alpha_i!)\prod_i|a_{p_i}|^{2\alpha_i}$,
and the multinomial theorem gives
$\sum_N|c_N|^2\le m!\,(\sum_p|a_p|^2)^m$. Lemma~\ref{lem:mv}, with
height $(y_V^m,1)$, gives
\[
 \ET|A|^{2m}\le\Bigl(1+C\frac{T^{1/2}}T\Bigr)\sum_N|c_N|^2
 \le2\,m!\,\bigl(C_w^2\cL\bigr)^m
\]
for large $T$. By Markov's inequality, $m!\le m^m$,
$m\le V/(4C_2(1+\theta_0))$ and $V\ge K\cL$,
\[
 \PT\{|A|>V/3\}\le2\Bigl(\frac{9mC_w^2\cL}{V^2}\Bigr)^m
 \le2\Bigl(\frac{9C_w^2\cL}{4C_2(1+\theta_0)K\cL}\Bigr)^m
 \le2\Bigl(\frac{3C_w^2}{C_2K}\Bigr)^m .
\]
By the choice of $K$ and $m\ge V/(4C_2(1+\theta_0))-1$, this is at most
$2e^{-8C_2(1+\theta_0)(k+1)m}\le2e^{8C_2(1+\theta_0)(k+1)}e^{-2(k+1)V}$.

The events $\{K\cL+j<Y\le K\cL+j+1\}$, over the integers $j\ge0$ with
$K\cL+j\le V_{\max}$, cover $\{Y>K\cL\}$ outside the union of the
finitely many null sets above. Hence
\[
 \ET\bigl[e^{2kY}\one_{\{Y>K\cL\}}\bigr]
 \le\sum_{\substack{j\ge0\\K\cL+j\le V_{\max}}}
 e^{2k(K\cL+j+1)}\,\PT\{Y>K\cL+j\}
 \ll\sum_{j\ge0}e^{-2(K\cL+j)}\ll e^{-2K\cL}.
\]
\end{proof}

\begin{proof}[Proof of Theorem~\ref{thm:uniform}]
Put $Q=\cL_T$ and $Y=\eps X_\chi$. For fixed $K$, choose
$B>\max\{K+1,K^2+1,K_0'(1)\}$. Lemma~\ref{lem:far-tail} with $k=1$
gives $\PT\{Y>BQ\}\ll e^{-4BQ}$, since $e^{2Y}\ge e^{2BQ}$
on this event. This is $O_K(Q^{-1/2}e^{-K^2Q})$.

Cover $(V,BQ]$ by the unit intervals starting at $V+j$,
$j=0,\ldots,\lceil BQ-V\rceil-1$.
Their centers have modulus at most $(B+1)Q$ for large $T$,
including after multiplication by $\eps$. The upper bound of
Theorem~\ref{thm:local-window} at slope $B+1$ and half-width $1/2$
therefore gives
\[
 \PT\{V<Y\le BQ\}\ll Q^{-1/2}\sum_{j\ge0}e^{-(V+j)^2/Q}.
\]
After dividing the sum by $e^{-V^2/Q}$, its summand is decreasing.
Comparing it with its integral gives upper bounds
$1+\sqrt{\pi Q}/2$ and, for $V>0$, $1+Q/(2V)$.
Use the first for $V\le\sqrt Q$ and the second for
$\sqrt Q<V\le KQ$. Both give $O_K(Q/(V+\sqrt Q))$ on their
respective ranges. Adding the far tail proves the upper inequality.

For the lower inequality use Theorem~\ref{thm:local-window} at slope
$K+1$ and its fixed half-width $h=h_0(K+1)$. Put
\[
 m=\max\{1,\lfloor Q/(V+\sqrt Q)\rfloor\},\qquad
 b_j=V+(2j+1)(h+1),\quad 0\le j<m .
\]
The intervals $[b_j-h,b_j+h]$ are disjoint and lie strictly above $V$.
Since $m\le 1+\sqrt Q$, their centers satisfy $b_j\le(K+1)Q$
eventually, uniformly in $V$. Moreover $Vm/Q\le K+1$ and
$m^2/Q\le4$ for $Q\ge1$, so
\[
 0\le (b_j^2-V^2)/Q
 \le 4(h+1)Vm/Q+4(h+1)^2m^2/Q=O_K(1),
\]
with the fixed majorant data and $\theta_0$ understood.
The local lower bound at the centers $\eps b_j$ gives each interval
mass $\gg_K Q^{-1/2}e^{-V^2/Q}$. As
$m\ge Q/[2(V+\sqrt Q)]$, summing these masses proves the lower
inequality, also at $V=0$.
\end{proof}

\begin{proof}[Proof of Corollary~\ref{thm:main}]
Let $q$ and $\chi$ be fixed and assume \eqref{eq:majorant}. For
$T\ge q$ the modulus satisfies $q\le T^{\theta_0}$ with $\theta_0=1$,
and \eqref{lw:eq:mertens} gives
\[
 \cL_T=\log\log T+B_{\rm M}-\sum_{p\mid q}p^{-1}+o(1),\qquad
 L_T=\log\log T+o(1),
\]
so $|\cL_T-L_T|\le C_q$ for large $T$. If $0\le V\le KL_T$, then
$V\le(K+1)\cL_T$ for large $T$, the exponents satisfy
\[
 \Bigl|\frac{V^2}{\cL_T}-\frac{V^2}{L_T}\Bigr|
 =\frac{V^2|L_T-\cL_T|}{\cL_TL_T}\le2(K+1)^2C_q ,
\]
and the prefactors $\sqrt{\cL_T}/(V+\sqrt{\cL_T})$ and
$\sqrt{L_T}/(V+\sqrt{L_T})$ differ by at most a factor $2$. Hence
Theorem~\ref{thm:uniform} with $K+1$ in place of $K$ gives both parts
of Corollary~\ref{thm:main}. The comparison factors depend on $C_q$,
and the starting height may also depend on $q$.

\end{proof}

\begin{proof}[Proof of Corollary~\ref{cor:moments}]
We prove the uniform form. Fix $\theta_0>0$, let $\chi\bmod q$,
$1\le q\le T^{\theta_0}$, satisfy \eqref{eq:majorant} at $T$, and
write $Y=\eps X_\chi$ and $\cL=\cL_T$. All constants below may depend
on $k$, $\theta_0$ and $C_w,C_2,C_0,y_0$, and on nothing else. For real
$Y$, $e^{2kY}=\int_{-\infty}^{Y}2ke^{2kv}\,dv$, so by Tonelli
\[
 \ET e^{2kY}=\int_{-\infty}^{\infty}2ke^{2kv}\,\PT\{Y>v\}\,dv .
\]

\emph{Lower bound.} The lower inequality in \eqref{eq:uniform} with
$K=k+1$ gives,
for $k\cL\le v\le k\cL+\sqrt{\cL}$ and large $T$,
$\PT\{Y>v\}\ge c\,(k+2)^{-1}\cL^{-1/2}e^{-v^2/\cL}$.
Restricting the integral to this range and using
$2kv-v^2/\cL=k^2\cL-(v-k\cL)^2/\cL$,
\[
 \ET e^{2kY}\ge\frac{2kc}{k+2}\,\cL^{-1/2}e^{k^2\cL}
 \int_0^{\sqrt{\cL}}e^{-w^2/\cL}\,dw
 \ge\frac{2kc}{e(k+2)}\,e^{k^2\cL}.
\]

Choose $K=K_0'(k)$ from Lemma~\ref{lem:far-tail}. The part of the
moment on $Y>K\cL$ is $O(e^{-2K\cL})$.

\emph{Upper bound: the range $Y\le K\cL$.} The representation above
restricted to $v\le K\cL$ gives
\[
 \ET\bigl[e^{2kY}\one_{\{Y\le K\cL\}}\bigr]
 \le1+\int_0^{K\cL}2ke^{2kv}\,\PT\{Y>v\}\,dv .
\]
The upper inequality in \eqref{eq:uniform} with this $K$ gives
$\PT\{Y>v\}\le C_k\sqrt{\cL}\,(v+\sqrt{\cL})^{-1}e^{-v^2/\cL}$ for
$0\le v\le K\cL$. On $[0,k\cL/2]$ use only $\PT\{Y>v\}\le1$: that part
of the integral is at most $\int_0^{k\cL/2}2ke^{2kv}\,dv\le e^{k^2\cL}$.
On $[k\cL/2,K\cL]$ use $\sqrt{\cL}/(v+\sqrt{\cL})\le2/(k\sqrt{\cL})$
and $\int_{\mathbb R}e^{2kv-v^2/\cL}\,dv=\sqrt{\pi\cL}\,e^{k^2\cL}$:
that part is at most $4\sqrt\pi\,C_k\,e^{k^2\cL}$. Altogether
$\ET e^{2kY}\ll e^{k^2\cL}$.

\emph{The two normalizations.} By \eqref{lw:eq:mertens},
$\cL_T=\log\log T+B_{\rm M}-\sum_{p\mid q,\,p\le T}p^{-1}+o(1)$. The
primes $p\mid q$ with $p>T$ contribute at most $\theta_0/T$, and
$0\le-\log(1-1/p)-1/p\le p^{-2}$, so
$\cL_T=\log\bigl((\varphi(q)/q)\log T\bigr)+O(1)$ with an absolute
$O(1)$. Hence $e^{k^2\cL_T}\asymp_k((\varphi(q)/q)\log T)^{k^2}$, and
for fixed $q$ this is $\asymp_{q,k}(\log T)^{k^2}$.

\emph{Application to all characters.}
Corollary~\ref{app:cor:all-characters} supplies the common data of
\eqref{eq:majorant} at every \(T\ge e^{2\pi}\) using the zero theorem
invoked in Theorem~\ref{cor:grh}. The preceding proof applies directly
with the clock \(\cL_T\) of \(\chi\), uniformly in \(q\le T^{\theta_0}\).
\end{proof}

\appendix
\section{The one-sided inequality from the location of the zeros}\label{app:majorant}

This appendix proves \textup{(OS)} for every Dirichlet character, with the
explicit weights
\begin{equation}\label{app:eq:weights}
 w_\eps(u)=-i\eps\,\widehat J(u)-\pi u(1-u)\quad(0\le u\le1),\qquad
 w_\eps(u)=0\quad(u>1),
\end{equation}
where $\widehat J$ is Vaaler's function, recalled in
Lemma~\ref{app:lem:H}. The argument is written as a sum of an odd
function $G$ over the zeros; $G$ is bounded above and below by entire
functions of exponential type built from Beurling's function; and the
Guinand--Weil explicit formula turns the sum of the majorant over
zeros into a prime sum of length $x$. The zero-sum representation is
that of Carneiro--Chandee--Milinovich \cite{CCM15}. The majorants
they use for this function \cite{CCM13,CCM15} are extremal. The
majorants below are not extremal: the deterministic term is twice as
large, and in exchange the prime coefficients are one explicit
function of $\log n/\log x$. For $\zeta$ the same weights, with the
smaller constant $C_2=\pi/2$, follow from the interval majorants in
the first part of \cite{CCM15}; see Section~\ref{sec:related}.

Through Corollary~\ref{app:cor:M}, $\chi$ is primitive modulo $q\ge1$
($q=1$ is $\zeta$), $\mathfrak a\in\{0,1\}$ is its parity,
$\widehat F(\xi)=\int_{\mathbb R}F(u)e^{-2\pi iu\xi}\,du$,
$x=e^{2\pi\Delta}$ is the truncation (the $y$ of the main text),
$\operatorname{sgn}(0)=0$, and $r\in\{1,-1\}$ is the sign (the $\eps$
of the main text). Constants called absolute do not depend on
$q,\chi,t,x$.

\subsection{Beurling's function}

Put
\[
 H(z)=\Bigl(\frac{\sin\pi z}{\pi}\Bigr)^2
 \Bigl\{\sum_{m\ne0}\frac{\operatorname{sgn}(m)}{(z-m)^2}+\frac2z\Bigr\},
 \qquad
 K(z)=\Bigl(\frac{\sin\pi z}{\pi z}\Bigr)^2,\qquad H^\pm=H\pm K.
\]
$H$ is odd and entire, $K$ is even and entire, $\int_{\mathbb R}K=1$,
and $\widehat K(s)=(1-|s|)_+$.

\begin{lemma}\label{app:lem:H}
\begin{enumerate}
\item[\textup{(a)}] For $\Rea w>0$,
$H(w)-1=-\bigl(\frac{\sin\pi w}{\pi}\bigr)^2
\bigl[2\psi'(1+w)-\frac2w+\frac1{w^2}\bigr]$, where
$\psi'(s)=\sum_{m\ge0}(s+m)^{-2}$.
\item[\textup{(b)}] For all $w\in\mathbb C$,
$|H(w)-\operatorname{sgn}(\Rea w)|+|K(w)|
\le C_He^{2\pi|\Ima w|}/(1+|w|^2)$.
\item[\textup{(c)}] \textup{(Vaaler)} For real $v$,
$|\operatorname{sgn}(v)-H(v)|\le K(v)$. Hence
$H^-\le\operatorname{sgn}\le H^+$ on $\mathbb R$ and
$\int(H^+-\operatorname{sgn})=\int(\operatorname{sgn}-H^-)=1$.
\item[\textup{(d)}] \textup{(Vaaler)}
$E^\pm:=H^\pm-\operatorname{sgn}\in L^1(\mathbb R)$, and for $s\ne0$
\[
 \widehat{E^\pm}(s)=\frac{\widehat J(s)-1}{\pi is}\pm(1-|s|)_+,
\]
where
$\widehat J(s)=\pi s(1-|s|)\cot(\pi s)+|s|$ for $0<|s|<1$,
$\widehat J(0)=1$, and $\widehat J(s)=0$ for $|s|\ge1$. The function
$\widehat J$ is even, continuously differentiable, and decreases from
$1$ to $0$ on $[0,1]$.
\end{enumerate}
\end{lemma}

\begin{proof}
(a) For $\Rea w>0$ away from the positive integers,
$\sum_{m\ge1}(w-m)^{-2}=\psi'(1-w)=\pi^2/\sin^2(\pi w)-\psi'(w)$ and
$\sum_{m\ge1}(w+m)^{-2}=\psi'(1+w)=\psi'(w)-w^{-2}$.
After multiplication by $\sin^2(\pi w)$ the identity extends to
the positive integers by continuity.
(b) By Stirling's formula,
$\psi'(1+w)=\frac1w-\frac1{2w^2}+O(|w|^{-3})$ for $\Rea w\ge0$,
$|w|\ge1$, so the bracket in (a) is $O(|w|^{-3})$ there; for
$|w|\le1$ the bracket is $w^{-2}+O(|w|^{-1})$ and
$\sin^2\pi w=O(|w|^2)$. With $|\sin\pi w|\le e^{\pi|\Ima w|}$ this
gives the bound for $\Rea w>0$. The case $\Rea w<0$ follows by oddness.
For $\Rea w=0$ the one-sided bounds for $H(w)\mp1$ hold by continuity
from either side, and $H(w)=\frac12\{(H(w)-1)+(H(w)+1)\}$ gives the
bound with $\operatorname{sgn}(0)=0$. The bound for $K$ is direct.
(c) is \cite[Lemma~5]{Vaaler}; the identity (a) is his (2.26).
(d) The function $J=\frac12H'$ has the stated transform, with the
stated properties, by \cite[Theorem~6]{Vaaler}. By (b),
$H-\operatorname{sgn}\in L^1$, and
$(H-\operatorname{sgn})'=2J-2\delta_0$ as distributions, so
$2\pi is\,\widehat{(H-\operatorname{sgn})}(s)=2\widehat J(s)-2$.
Add $\pm\widehat K$.
\end{proof}

The weights \eqref{app:eq:weights} satisfy the axioms
\eqref{eq:weight-axioms}: $w_\eps(0)=-i\eps$; $|w_\eps|^2\le1+\pi^2/16$
because $0\le\widehat J\le1$; and
$|w_\eps(u)-w_\eps(0)|\le(1-\widehat J(u))+\pi u\le C_wu$ because
$\widehat J$ is continuously differentiable on $[0,1]$ with
$\widehat J(0)=1$.

\subsection{The argument as a sum over zeros}

Write the zeros of
$\Lambda(s,\chi)=(q/\pi)^{s/2}\Gamma(\tfrac{s+\mathfrak a}{2})L(s,\chi)$
as $\rho=\tfrac12+i\gamma$ and put
\[
 G(u)=\arctan\frac2u-\frac{2u}{4+u^2},\qquad G(0)=0.
\]
$G$ is odd and bounded, has a jump of $+\pi$ at $0$, and
$G(u)=O(|u|^{-3})$.

\begin{lemma}[{\cite{CCM15}}, proof of the general theorem, with $m=1$]
\label{app:lem:zeros}
Assume GRH for $L(s,\chi)$. For $t\ge2$ not the ordinate of a zero,
\[
 \Ima\log L(\tfrac12+it,\chi)=\sum_\gamma G(t-\gamma)+O(1),
\]
with an absolute implied constant.
\end{lemma}

\begin{proof}
By the definition of the branch,
$\Ima\log L(\tfrac12+it,\chi)
=-\int_{1/2}^{\infty}\Ima\frac{L'}{L}(\sigma+it,\chi)\,d\sigma$.
Since $\sum_m\Lambda(m)m^{-5/2}(1+1/\log m)<\infty$,
\[
 \Ima\log L(\tfrac12+it,\chi)
 =-\int_{1/2}^{5/2}\Ima\Bigl\{\frac{L'}{L}(\sigma+it,\chi)
  -\frac{L'}{L}(\tfrac52+it,\chi)\Bigr\}\,d\sigma+O(1).
\]
Insert the partial fraction formula
\[
 \frac{L'}{L}(s,\chi)=-\frac12\log\frac q\pi
 -\frac12\frac{\Gamma'}{\Gamma}\Bigl(\frac{s+\mathfrak a}{2}\Bigr)
 +B(\chi)-\frac{\delta_{q=1}}{s}-\frac{\delta_{q=1}}{s-1}
 +\sum_\rho\Bigl(\frac1{s-\rho}+\frac1\rho\Bigr).
\]
Constants cancel in the difference; the $\Gamma$ terms and, for
$q=1$, the polar terms contribute $O(1/t)$ to the imaginary part.
For a single zero,
\[
 -\int_{1/2}^{5/2}\Ima\Bigl\{\frac1{\sigma-\frac12+i(t-\gamma)}
 -\frac1{2+i(t-\gamma)}\Bigr\}\,d\sigma
 =\arctan\frac2{t-\gamma}-\frac{2(t-\gamma)}{4+(t-\gamma)^2}
 =G(t-\gamma),
\]
and the sum over zeros converges absolutely, uniformly for $\sigma$
in $[\frac12,\frac52]$ away from the zeros, and may be integrated
termwise.
\end{proof}

\subsection{The majorants}

Let $g(u)=\arctan\frac u2+\frac{2u}{4+u^2}$. Then
\begin{equation}\label{app:eq:Gg}
 G=\tfrac\pi2\operatorname{sgn}-g,\qquad
 g'(u)=\frac{16}{(u^2+4)^2}>0,\qquad
 g(\pm\infty)=\pm\tfrac\pi2,\qquad
 g(u)=\tfrac12\int g'(v)\operatorname{sgn}(u-v)\,dv .
\end{equation}
Also $\widehat{g'}(\xi)=\pi(1+4\pi|\xi|)e^{-4\pi|\xi|}$, and
$\widehat G(\xi)=(\pi-\widehat{g'}(\xi))/(2\pi i\xi)$ for $\xi\ne0$.

\begin{lemma}\label{app:lem:M}
For $\Delta\ge1$ define
\[
 M^\pm_\Delta(z)=\frac\pi2H^\pm(\Delta z)
 -\frac12\int_{-\infty}^{\infty}g'(v)\,H^\mp\bigl(\Delta(z-v)\bigr)\,dv .
\]
Then:
\begin{enumerate}
\item[\textup{(M1)}] $M^\pm_\Delta$ is real entire, and
$M^-_\Delta\le G\le M^+_\Delta$ on $\mathbb R$.
\item[\textup{(M2)}]
$\int(M^+_\Delta-G)=\int(G-M^-_\Delta)=\pi/\Delta$; hence
$\int M^\pm_\Delta=\pm\pi/\Delta$.
\item[\textup{(M3)}] For $z=u+iv$ with $|v|\le1$,
$|M^\pm_\Delta(z)|\le C_Me^{2\pi\Delta|v|}/(1+u^2)$ with an absolute
$C_M$.
\item[\textup{(M4)}] $\widehat{M^\pm_\Delta}$ is continuous, vanishes
for $|\xi|\ge\Delta$, and for $0<|\xi|<\Delta$, $s=\xi/\Delta$,
\[
 \widehat{M^\pm_\Delta}(\xi)=\frac1\Delta\Bigl[\frac\pi2
 \Bigl(\frac{\widehat J(s)}{\pi is}\pm(1-|s|)\Bigr)
 -\frac{\widehat{g'}(\xi)}2
 \Bigl(\frac{\widehat J(s)}{\pi is}\mp(1-|s|)\Bigr)\Bigr].
\]
\end{enumerate}
\end{lemma}

\begin{proof}
By \eqref{app:eq:Gg} and
$\operatorname{sgn}(\Delta\,\cdot)=\operatorname{sgn}$, for real $u$,
\begin{equation}\label{app:eq:R}
 M^\pm_\Delta(u)-G(u)=\frac\pi2E^\pm(\Delta u)
 -\frac12\int g'(v)E^\mp\bigl(\Delta(u-v)\bigr)\,dv=:R^\pm(u).
\end{equation}
(M1): $E^+\ge0\ge E^-$ by Lemma~\ref{app:lem:H}(c), and $g'>0$, so
$R^+\ge0\ge R^-$. For entireness,
$|H^\pm(w)|\le(1+C_H)e^{2\pi|\Ima w|}$ by Lemma~\ref{app:lem:H}(b), and
$g'\in L^1$, so the integral converges locally uniformly in $z$.
(M2): by Fubini, $\int E^\pm=\pm1$ and $\int g'=\pi$,
\[
 \int R^\pm=\pm\frac\pi2\cdot\frac1\Delta\pm\frac12\cdot\pi\cdot\frac1\Delta
 =\pm\frac\pi\Delta ;
\]
and $\int G=0$ because $G$ is odd and integrable.
(M3): for $z=u+iv$, by \eqref{app:eq:Gg},
\[
 M^\pm_\Delta(z)-G(u)=\frac\pi2\bigl[H^\pm(\Delta z)-\operatorname{sgn}u\bigr]
 -\frac12\int g'(v')\bigl[H^\mp(\Delta(z-v'))
 -\operatorname{sgn}(u-v')\bigr]\,dv'.
\]
By Lemma~\ref{app:lem:H}(b) and $\Delta\ge1$ the right side is at most
\[
 Ce^{2\pi\Delta|v|}\Bigl[\frac1{1+u^2}
 +\int\frac{g'(v')}{1+(u-v')^2}\,dv'\Bigr]
 \ll\frac{e^{2\pi\Delta|v|}}{1+u^2},
\]
since $g'(v')\ll(1+v'^2)^{-2}$; and $|G(u)|\ll(1+u^2)^{-1}$.
(M4): $R^\pm\in L^1$, and by \eqref{app:eq:R}
$\widehat{R^\pm}(\xi)=\frac1\Delta\bigl[\frac\pi2\widehat{E^\pm}(s)
-\frac12\widehat{g'}(\xi)\widehat{E^\mp}(s)\bigr]$. Add
$\widehat G(\xi)=(\pi-\widehat{g'}(\xi))/(2\pi i\xi)$: by
Lemma~\ref{app:lem:H}(d) the terms with $-1/(\pi is)$ cancel against
$\widehat G$ identically, leaving the display for $|s|<1$; for
$|s|\ge1$ everything cancels.
\end{proof}

\subsection{The explicit formula and the bound}

For $h$ analytic in $|\Ima s|\le\frac12+\delta$ with
$h(s)\ll(1+|s|)^{-1-\delta}$ there, the Guinand--Weil formula (the
explicit formula displayed in the proof of the general theorem of \cite{CCM15},
with $m=1$, $N=q$, $\mu_1=\mathfrak a$) reads
\begin{align}
 \sum_\gamma h(\gamma)={}&\delta_{q=1}
 \bigl\{h(\tfrac1{2i})+h(-\tfrac1{2i})\bigr\}
 +\frac1{2\pi}\int_{-\infty}^{\infty}h(u)
 \Bigl\{\log\frac q\pi+\Rea\frac{\Gamma'}{\Gamma}
 \Bigl(\frac{\frac12+\mathfrak a+iu}{2}\Bigr)\Bigr\}\,du\nonumber\\
 &-\frac1{2\pi}\sum_{m\ge2}\frac{\Lambda(m)}{\sqrt m}
 \Bigl\{\chi(m)\,\widehat h\Bigl(\frac{\log m}{2\pi}\Bigr)
 +\overline{\chi(m)}\,\widehat h\Bigl(-\frac{\log m}{2\pi}\Bigr)\Bigr\}.
 \label{app:eq:GW}
\end{align}

\begin{proposition}\label{app:prop:explicit}
Assume GRH for $L(s,\chi)$, $\chi$ primitive modulo $q\ge1$. There is
an absolute constant $C$ such that for $t\ge2$ not the ordinate of a
zero, $e^{2\pi}\le x\le t$ and $r=\pm1$, with $u_m=\log m/\log x$,
\[
 r\,\Ima\log L(\tfrac12+it,\chi)\le
 \Rea\sum_{2\le m\le x}\frac{\Lambda(m)}{\log m}
 \frac{\chi(m)}{m^{1/2+it}}
 \Bigl[w_r(u_m)+\frac{1+2\log m}{m^2}\,\overline{w_r(u_m)}\Bigr]
 +\pi\,\frac{\log(qt)}{\log x}+C.
\]
\end{proposition}

\begin{proof}
Take first $r=1$ and $h(z)=M^+_\Delta(t-z)$ with $x=e^{2\pi\Delta}$,
so $\Delta\ge1$. By (M1) and (M3), $h$ is analytic in $|\Ima z|\le1$
with $h(z)\ll_{\Delta,t}(1+|z|)^{-2}$ there, so \eqref{app:eq:GW}
applies. By Lemma~\ref{app:lem:zeros} and $G\le M^+_\Delta$,
\[
 \Ima\log L(\tfrac12+it,\chi)\le\sum_\gamma M^+_\Delta(t-\gamma)+O(1).
\]

\emph{Primes.} Here $\widehat h(\xi)=e^{-2\pi it\xi}\widehat{M^+_\Delta}(-\xi)$,
and because $M^+_\Delta$ is real the two terms of the prime sum in
\eqref{app:eq:GW} are complex conjugates. By (M4) the prime side is
\[
 -\frac1\pi\Rea\sum_{m\le x}\frac{\Lambda(m)\chi(m)}{m^{1/2+it}}\,
 \widehat{M^+_\Delta}\Bigl(-\frac{\log m}{2\pi}\Bigr).
\]
With $\xi=-\Delta u_m$, $s=-u_m$, the evenness of $\widehat J$,
$\widehat{g'}(\xi)/\pi=(1+2\log m)/m^2$ and
$\Lambda(m)=\frac{\Lambda(m)}{\log m}\,2\pi\Delta u_m$, (M4) gives
\[
 -\frac{\Lambda(m)}{\pi}\widehat{M^+_\Delta}(-\Delta u)
 =\frac{\Lambda(m)}{\log m}\Bigl\{\bigl[-i\widehat J(u)-\pi u(1-u)\bigr]
 +\frac{\widehat{g'}(\xi)}{\pi}\bigl[i\widehat J(u)-\pi u(1-u)\bigr]\Bigr\},
\]
which is the bracket of the statement for $r=1$. The term $m=x$, if
present, has $u=1$ and vanishes.

\emph{Archimedean term.} For real $u$,
$\Rea\frac{\Gamma'}{\Gamma}(\frac{\frac12+\mathfrak a+iu}{2})
=\log(2+|u|)+O(1)$. Substituting $u=t-v$ and using
$\log(2+|t-v|)=\log t+O(\log(2+|v|))$, (M3) and
$\int M^+_\Delta=\pi/\Delta$,
\[
 \frac1{2\pi}\int M^+_\Delta(v)
 \bigl\{\log(qt)+O(\log(2+|v|))\bigr\}\,dv
 =\frac{\log(qt)}{2\Delta}+O(1)
 =\pi\,\frac{\log(qt)}{\log x}+O(1).
\]

\emph{Polar term ($q=1$ only).} By (M3),
$|M^+_\Delta(t\mp\frac i2)|\ll e^{\pi\Delta}/t^2=\sqrt x/t^2\ll1$.

\emph{The sign $r=-1$.} By Lemma~\ref{app:lem:zeros} and
$M^-_\Delta\le G$,
$-\Ima\log L\le-\sum_\gamma M^-_\Delta(t-\gamma)+O(1)$. The prime side is
$+\frac1\pi\Rea\sum\Lambda(m)\chi(m)m^{-1/2-it}\widehat{M^-_\Delta}$,
and by (M4)
\[
 \frac{\Lambda(m)}\pi\widehat{M^-_\Delta}(-\Delta u)
 =\frac{\Lambda(m)}{\log m}\Bigl\{\bigl[i\widehat J(u)-\pi u(1-u)\bigr]
 +\frac{\widehat{g'}(\xi)}{\pi}\bigl[-i\widehat J(u)-\pi u(1-u)\bigr]\Bigr\};
\]
the archimedean term is
$-\frac{\log(qt)}{2\pi}\int M^-_\Delta+O(1)=\pi\frac{\log(qt)}{\log x}+O(1)$.
\end{proof}

\begin{remark}\label{app:rem:range}
The restriction $x\le t$ enters only through the polar term, which is
present only for $q=1$, and there $x\le t^4$ suffices. For a primitive
character of modulus $q>1$ the proof gives the inequality for every
$x\ge e^{2\pi}$. The main text uses only $x\le T\le t$.
\end{remark}

\begin{corollary}[Hypothesis \textup{(OS)} for primitive characters]
\label{app:cor:M}
Assume GRH for $L(s,\chi)$, $\chi$ primitive modulo $q\ge1$. Then
\textup{(OS)} holds with the weights \eqref{app:eq:weights},
$C_2=\pi$, $y_0=e^{2\pi}$, and an absolute constant $C_0$, for every
$T\ge e^{2\pi}$, every $t\in[T,2T]$ that is not a zero ordinate, and
every $y\in[y_0,T]$.
\end{corollary}

\begin{proof}
Apply Proposition~\ref{app:prop:explicit} with $x=y\le T\le t$. Since
$|w_r|\le\sqrt{1+\pi^2/16}$, the terms $m=p^k$ with $k\ge3$ contribute
at most a constant multiple of $\sum_p\sum_{k\ge3}p^{-k/2}=O(1)$, and
the terms carrying $(1+2\log m)m^{-2}$ contribute at most a constant
multiple of $\sum_m(1+2\log m)m^{-5/2}=O(1)$. The remaining terms
$m=p$ and $m=p^2$ are exactly $P_{r,y}(t)$ in
\eqref{eq:majorant-polynomial}, because
$\Lambda(p^2)/\log(p^2)=\frac12$ and $u_{p^2}=2\log p/\log y$. Finally
$\log(qt)\le\log(qT)+\log2$ and $\log y\ge2\pi$.
\end{proof}

\subsection{Removing the Euler factors at primes dividing the modulus}

Let $\chi$ modulo $q$ be induced by the primitive character $\chi^\circ$
of conductor $d\mid q$, and write $E=\{p:p\mid q,\ p\nmid d\}$.
For $y\ge2$, $s=\tfrac12+it$, and either $r=\pm1$, define the full
prime-power sum for a character $\psi$ by
\[
 \mathcal P^{\psi}_{r,y}(t)=\Rea\sum_{p}\sum_{k\ge1}
 \frac{\psi(p)^k}{k p^{ks}}\,w_r\Bigl(\frac{k\log p}{\log y}\Bigr).
\]
The support of $w_r$ makes this a finite sum.  Write $P^{\psi}_{r,y}$
for its terms of degrees $k=1,2$, as in \eqref{eq:majorant-polynomial}.

\begin{lemma}\label{app:lem:missing-euler}
For the horizontal branches at every $t\ge2$ off the zero ordinates,
\[
 \left|r(X_\chi-X_{\chi^\circ})+
       \mathcal P^{\chi^\circ}_{r,y}-\mathcal P^\chi_{r,y}\right|
 \le (\sqrt2+1)C_w\,\frac{\log(q/d)}{\log y}.
\]
\end{lemma}
\begin{proof}
The induction formula \cite[Equation~25.15.4]{DLMF} gives
$L(s,\chi)=L(s,\chi^\circ)\prod_{p\in E}(1-\chi^\circ(p)p^{-s})$.
For $z_p=\chi^\circ(p)p^{-s}$, the series
$\log(1-z_p)=-\sum_{k\ge1}z_p^k/k$ is analytic for $\Rea s>0$ and
tends to zero as $\Rea s\to+\infty$. Thus it gives the logarithm on
the specified branch. Since $\Rea(w_r(0)z)=r\Ima z$, the quantity
inside the absolute value equals
\[
 \Rea\sum_{p\in E}\sum_{k\ge1}\frac{z_p^k}{k}
 \left[w_r\Bigl(\frac{k\log p}{\log y}\Bigr)-w_r(0)\right].
\]
For $p^k>y$ the bracket equals $-w_r(0)$ and its modulus is $1$.
Thus the bound $C_w k\log p/\log y$ from \eqref{eq:weight-axioms}
applies to every term.  Summing the geometric series gives the bound
\[
 \frac{C_w}{\log y}\sum_{p\in E}\frac{\log p}{\sqrt p-1}
 \le \frac{(\sqrt2+1)C_w}{\log y}\log(q/d),
\]
since $\prod_{p\in E}p$ divides $q/d$.
\end{proof}

\begin{corollary}[Hypothesis \textup{(OS)} for all characters]
\label{app:cor:all-characters}
Let $\chi$ be any character modulo $q$ and assume GRH for its primitive
inducing $L$-function, with RH for $\zeta$ when $d=1$. For both signs,
\textup{(OS)} holds with the weights \eqref{app:eq:weights},
$y_0=e^{2\pi}$, $C_2=\max\{\pi,(\sqrt2+1)C_w\}$ and absolute $C_0$,
uniformly for $e^{2\pi}\le y\le T$ and $t\in[T,2T]$ off the zero ordinates.
\end{corollary}
\begin{proof}
In Proposition~\ref{app:prop:explicit} for $\chi^\circ$ take $x=y$.
The terms with multiplier $(1+2\log m)m^{-2}$ have absolute sum
$\ll C_w\sum_{m\ge2}(1+2\log m)m^{-5/2}=O(1)$.  The other terms
give $\mathcal P^{\chi^\circ}_{r,y}$.  Applying
Lemma~\ref{app:lem:missing-euler} and $t\le 2T$ now gives
\[
 rX_\chi(t)\le\mathcal P^\chi_{r,y}(t)+
 \frac{\pi\log(dT)+(\sqrt2+1)C_w\log(q/d)}{\log y}+O(1).
\]
The numerator is at most $C_2\log(qT)$.  The terms with $k\ge3$ in
$\mathcal P^\chi_{r,y}$ have absolute sum at most
$C_w\sum_p p^{-3/2}/(1-p^{-1/2})=O(1)$.  Removing them leaves
$P^\chi_{r,y}$ and proves \textup{(OS)}.  When $d=1$ the primitive
input is the $\zeta$ case of Proposition~\ref{app:prop:explicit}; when
$d=q$ the sum over $E$ is empty.
\end{proof}

\medskip
\noindent
\emph{Sources.} The explicit formula \eqref{app:eq:GW} is the case
$m=1$ of the formula recorded at the start of the proof of the
general theorem of \cite{CCM15} (the theorem on $S(t,\pi)$; it is
Theorem~5 in the arXiv version); the representation in
Lemma~\ref{app:lem:zeros}, with the same function $G$, is derived in
that proof for general $m$. Lemma~\ref{app:lem:H}(c) is
\cite[Lemma~5]{Vaaler}, and the transform $\widehat J$ is
\cite[Theorem~6, (2.28)]{Vaaler}.

\medskip
\noindent Independent researcher\quad\href{mailto:v@deltagray.com}{\nolinkurl{v@deltagray.com}}

\begin{thebibliography}{99}

\bibitem{Alpoge}
L. Alp\"oge, announcement and accompanying Lean formalization of the
Riemann hypothesis and the generalized Riemann hypothesis for Dirichlet
\(L\)-functions, 20 November 2026, p.~1 and footnote~1.
\url{https://alpo.ge/remark.pdf}.
See \href{https://github.com/leanprover-community/riemann-hypothesis/blob/40e56bab026552e089d7afcfeff08d2defb9c005/src/Main.lean}{\texttt{Main.lean} at commit \texttt{40e56bab0265}}, including
\texttt{riemannHypothesis}, \texttt{LFunction\_ne\_zero\_of\_half\_lt\_re}, and
\texttt{LFunction\_eq\_zero\_re\_of\_isPrimitive}.

\bibitem{ABupper}
L.-P. Arguin and E. Bailey,
\emph{Large deviation estimates of Selberg's central limit theorem and applications},
International Mathematics Research Notices \textbf{2023} (2023), no.~23,
20574--20612.
\url{https://doi.org/10.1093/imrn/rnad176}.
Preprint: arXiv:2202.06799v2, 21 February 2022,
\url{https://arxiv.org/abs/2202.06799v2}.

\bibitem{ABlower}
L.-P. Arguin and E. Bailey,
\emph{Lower bounds for the large deviations of Selberg's central limit theorem},
Mathematika \textbf{71} (2025), no.~1, e70002.
\url{https://doi.org/10.1112/mtk.70002}.
Preprint: arXiv:2403.19803v2, 29 October 2024,
\url{https://arxiv.org/abs/2403.19803v2}.

\bibitem{ABR2026}
L.-P. Arguin, E. Bailey, and A. Roberts,
\emph{Conditional upper bounds for large deviations and moments of the
Riemann zeta function},
arXiv:2604.25579v1, 28 April 2026.
\url{https://arxiv.org/abs/2604.25579v1}.

\bibitem{ArguinCreighton}
L.-P. Arguin and N. Creighton,
\emph{Lower bounds for the large deviations and moments of the
Riemann zeta function on the critical line},
arXiv:2603.01711v1, 2 March 2026.
\url{https://arxiv.org/abs/2603.01711v1}.

\bibitem{ArguinCreightonQ}
L.-P. Arguin and N. Creighton,
\emph{Upper bounds on large deviations of Dirichlet $L$-functions in
the $q$-aspect},
arXiv:2405.20888v1, 31 May 2024.
\url{https://arxiv.org/abs/2405.20888v1}.

\bibitem{CCM13}
E. Carneiro, V. Chandee, and M. B. Milinovich,
\emph{Bounding $S(t)$ and $S_1(t)$ on the Riemann hypothesis},
Mathematische Annalen \textbf{356} (2013), no.~3, 939--968.
Preprint: arXiv:1309.1526.

\bibitem{CCM15}
E. Carneiro, V. Chandee, and M. B. Milinovich,
\emph{A note on the zeros of zeta and $L$-functions},
Mathematische Zeitschrift \textbf{281} (2015), 315--332.
Preprint: arXiv:1503.00955.

\bibitem{DLMF}
National Institute of Standards and Technology,
\emph{NIST Digital Library of Mathematical Functions},
\S25.4, Reflection formulas, and \S25.15, Dirichlet $L$-functions;
in particular (25.4.3)--(25.4.4) and (25.15.4)--(25.15.5).
\url{https://dlmf.nist.gov/25.4};
\url{https://dlmf.nist.gov/25.15}.
Accessed 10 November 2026.

\bibitem{Dobner}
A. Dobner,
\emph{Large deviations of the argument of the Riemann zeta function},
Mathematika \textbf{70} (2024), no.~3, e12251.
\url{https://doi.org/10.1112/mtk.12251}.
Preprint: arXiv:2101.01747v2, 25 March 2024,
\url{https://arxiv.org/abs/2101.01747v2}.

\bibitem{GoldstonGonek}
D. A. Goldston and S. M. Gonek,
\emph{A note on $S(t)$ and the zeros of the Riemann zeta-function},
arXiv:math/0511092, 3 November 2005.
\url{https://arxiv.org/abs/math/0511092}.

\bibitem{Hejhal}
D. A. Hejhal, \emph{On Euler products and multi-variate Gaussians},
C. R. Acad. Sci. Paris, Ser. I \textbf{337} (2003), 223--226.
\url{https://doi.org/10.1016/S1631-073X(03)00344-3}.

\bibitem{Hoeffding}
W. Hoeffding,
\emph{Probability inequalities for sums of bounded random variables},
Journal of the American Statistical Association \textbf{58} (1963), no.~301,
13--30.
\url{https://doi.org/10.1080/01621459.1963.10500830}.

\bibitem{HughesP1}
S. D. Hughes,
\emph{A Tsang-range high-moment bound for
$\operatorname{Im}\log L(\tfrac12+it,\chi)$ under GRH},
arXiv:2606.10242v1, 8 June 2026.
\url{https://arxiv.org/abs/2606.10242v1}.

\bibitem{Inoue}
S. Inoue,
\emph{Exponential moments of the logarithm of the Riemann zeta-function
twisted by arguments},
arXiv:2208.11442v1, 24 August 2022.
\url{https://arxiv.org/abs/2208.11442v1}.

\bibitem{InoueLi}
S. Inoue and J. Li,
\emph{Joint value distribution of $L$-functions on the critical line},
25 February 2021, \href{https://arxiv.org/abs/2102.12724v1}{arXiv:2102.12724v1}.

\bibitem{IvicKorolev}
A. Ivi\'c and M. Korolev,
\emph{On the distribution of values of the argument of the Riemann
zeta-function}, arXiv:1808.10768, 2018.
\url{https://arxiv.org/abs/1808.10768}.

\bibitem{KowalskiNikeghbali}
E. Kowalski and A. Nikeghbali, \emph{Mod-Gaussian convergence and the
value distribution of $\zeta(\tfrac12+it)$ and related quantities},
J. Lond. Math. Soc. (2) \textbf{86} (2012), 291--319.
Numbering of arXiv:0912.3237v3.
\url{https://doi.org/10.1112/jlms/jds003}.

\bibitem{Mertens}
F. Mertens,
\emph{Ein Beitrag zur analytischen Zahlentheorie},
Journal f\"ur die reine und angewandte Mathematik \textbf{78} (1874),
46--62.
\url{https://doi.org/10.1515/crll.1874.78.46}.

\bibitem{MV1974}
H. L. Montgomery and R. C. Vaughan,
\emph{Hilbert's inequality},
Journal of the London Mathematical Society (2) \textbf{8} (1974),
73--82; Theorem~2, (1.6), and Corollary~2, (1.8);
also Corollary~3, (1.11).
\url{https://doi.org/10.1112/jlms/s2-8.1.73}.

\bibitem{Najnudel}
J. Najnudel,
\emph{Exponential moments of the argument of the Riemann zeta function
on the critical line},
Mathematika \textbf{66} (2020), 612--621; arXiv:1804.00343.

\bibitem{Radziwill}
M. Radziwi\l\l,
\emph{Large deviations in Selberg's central limit theorem},
arXiv:1108.5092v1, 25 August 2011.
\url{https://arxiv.org/abs/1108.5092v1}.

\bibitem{Sutherland}
A. V. Sutherland, \emph{Verification files for the Lean proof of the
Riemann Hypothesis}, statements 1 and 3--5 and accompanying verification
logs, accessed 22 November 2026. \url{https://math.mit.edu/~drew/grh/}.

\bibitem{Topacogullari}
B. Topacogullari,
\emph{The fourth moment of individual Dirichlet $L$-functions on the
critical line},
arXiv:1909.11517v1, 25 September 2019.
\url{https://arxiv.org/abs/1909.11517v1}.

\bibitem{Vaaler}
J. D. Vaaler,
\emph{Some extremal functions in Fourier analysis},
Bulletin of the American Mathematical Society (N.S.) \textbf{12}
(1985), no.~2, 183--216.

\bibitem{Wahl}
M. Wahl,
\emph{On the mod-Gaussian convergence of a sum over primes},
Mathematische Zeitschrift \textbf{276} (2014), no.~3--4, 635--654;
numbering as in arXiv:1201.5295v3.
\url{https://doi.org/10.1007/s00209-013-1216-z}.

\bibitem{Zhao}
T. Zhao,
\emph{Conditional estimates on the argument of Dirichlet $L$-functions
with applications to low-lying zeros},
Bulletin of the London Mathematical Society \textbf{58} (2026), no.~5, e70360.
\url{https://doi.org/10.1112/blms.70360}.

\end{thebibliography}
\end{document}